\documentclass{article}
\usepackage{amssymb,amsmath,amsthm}
\usepackage{fancyhdr}
\usepackage[british]{babel}
\usepackage{color}
\usepackage[title]{appendix}
\usepackage{amsmath,amsthm}

\usepackage{enumerate,enumitem}
\usepackage{enumitem,verbatim}

\usepackage{ulem}
\usepackage{hyperref}
\usepackage[margin=2.5cm]{geometry}
\usepackage{lineno}
\usepackage[ruled]{algorithm2e}
\usepackage{todonotes}

\newtheorem{theorem}{Theorem}[section]
\newtheorem{prop}[theorem]{Proposition}
\newtheorem{lemma}[theorem]{Lemma}
\newtheorem{cor}[theorem]{Corollary}

\theoremstyle{definition}

\newtheorem{question}{Question}

\newtheorem{defn}[theorem]{Definition}

\newtheorem{claim}[theorem]{Claim}
\newtheorem{fact}[theorem]{Fact}

\newcommand{\btheorem}{\begin{theorem}}
\newcommand{\etheorem}{\end{theorem}}
\newcommand{\bconjecture}{\begin{conjecture}}
\newcommand{\econjecture}{\end{conjecture}}
\newcommand{\bproposition}{\begin{proposition}}
\newcommand{\eproposition}{\end{proposition}}
\newcommand{\bdefinition}{\begin{definition}}
\newcommand{\edefinition}{\end{definition}}
\newcommand{\bcorollary}{\begin{corollary}}
\newcommand{\ecorollary}{\end{corollary}}
\newcommand{\bproof}{\begin{proof}}
\newcommand{\eproof}{\end{proof}}
\newcommand{\bclaim}{\begin{claim}}
\newcommand{\eclaim}{\end{claim}}
\newcommand{\bquestion}{\begin{question}}
\newcommand{\equestion}{\end{question}}
\newcommand{\bfact}{\begin{fact}}
\newcommand{\efact}{\end{fact}}
\newcommand{\bremark}{\begin{remark}}
\newcommand{\eremark}{\end{remark}}
\newcommand{\eexample}{\end{example}}
\newcommand{\bexample}{\begin{example}}

\newcommand{\elemma}{\end{lemma}}
\newcommand{\blemma}{\begin{lemma}}

\newcommand{\N}{\mathbb{N}}

\newcommand{\eps}{\varepsilon}
\newcommand{\<}{\subseteq}

\newcommand{\red}[1]{\textcolor[rgb]{1.00,0.00,0.00}{#1}}
\newcommand{\De}{\Delta}
\newcommand{\de}{\delta}

\title{Transversal packings in families of percolated hypergraphs}

\author{Jie Han\thanks{School of Mathematics and Statistics and Center for Applied Mathematics, Beijing Institute of Technology, Beijing, China. Email: {\tt han.jie@bit.edu.cn}.},
Jie Hu\thanks{School of Mathematical Sciences and LPMC, Nankai University, Tianjin, China. Email: {\tt hujie@nankai.edu.cn}.},
Shunan Wei\thanks{School of Mathematics, Shandong University, Jinan, China. Email: {\tt snwei@mail.sdu.edu.cn}.},
Donglei Yang\thanks{School of Mathematics, Shandong University, Jinan, China. Email: {\tt dlyang@sdu.edu.cn}.}}

\date{\today}

\begin{document}
\maketitle
\linespread{1.2}

\begin{abstract}
  Let $F$ be a strictly $1$-balanced $k$-graph on $s$ vertices with $t$ edges and $\delta_{F,d}^T$ be the infimum of $\delta>0$ such that for every $\alpha>0$ and sufficiently large $n\in \mathbb{N}$, every $k$-graph system $\mathbf H=\{H_{1}, H_{2}, \dots ,H_{tn}\}$ on the same $sn$ vertices with $\delta_d(H_i)\ge (\delta+\alpha)\binom{sn-d}{k-d}$, $i\in [tn]$ contains a transversal $F$-factor, that is, an $F$-factor consisting of exactly one edge from each $H_i$.
  In this paper we prove the following result. Let $\mathbf{H} =\{H_{1}, H_{2}, \dots ,H_{tn}\}$ be a   $k$-graph system where each $H_{i}$ is an $sn$-vertex $k$-graph with $\delta_d(H_i)\ge (\delta_{F,d}^T+\alpha)\binom{sn-d}{k-d}$. Then with high probability $\mathbf{H}(p) :=\{H_{1}(p), H_{2}(p), \dots ,H_{tn}(p)\}$ contains a transversal $F$-factor, where  $H_i(p)$ is a random subhypergraph of $H_i$ and 
  $p=\Omega(n^{-1/d_1(F)-1}(\log n)^{1/t})$.
  This extends a recent result by Kelly, M\"{u}yesser and Pokrovskiy, and independently by Joos, Lang and Sanhueza-Matamala. Moreover, the assumption on $p$ is best possible up to a constant.
  Along the way, we also obtain a spread version of a result of Pikhurko on perfect matchings in $k$-partite $k$-graphs.
  
\end{abstract}


\section{Introduction}
Determining the minimum degree condition for the existence of spanning structures is a central
theme in extremal graph theory, which is called the \textit{Dirac-type} problem because of the cornerstone result of Dirac~\cite{Dirac} on the existence of Hamilton cycles.
Given graphs $H$ and $G$, an \textit{$H$-tiling} of $G$ is a collection of vertex-disjoint copies of $H$ in $G$, which naturally extends the definition of matchings. An \textit{$H$-factor} of $G$ is an $H$-tiling which covers all vertices of $G$.
Classical results of Corr\'adi and Hajnal~\cite{Corradi}, and Hajnal and Szemer\'edi~\cite{HSz} determined the minimum degree for the existence of a $K_r$-factor.

\begin{theorem}[\cite{Corradi,HSz}]
\label{thm:HSz}
Let $G$ be an $n$-vertex graph with $n\in r\mathbb{N}$. If $\delta(G)\ge \left(1-1/r\right)n$, then $G$ contains a $K_r$-factor.
\end{theorem}

Since then there has been a significant body of research on determining the minimum degree threshold forcing $F$-factors in graphs, hypergraphs, directed graphs, etc.
The problem for graphs has been settled almost completely by K\"uhn and Osthus~\cite{MR2506388}, who obtained,
up to an additive constant, the best possible minimum degree condition which forces
an $F$-factor.

\subsection{Robustness}
The study of the \textit{robustness} of graph properties has received considerable attention recently, aiming to enhance classic results in extremal graph theory and probabilistic combinatorics.
We refer the reader to a comprehensive survey of Sudakov~\cite{Sudakov-survey} which collects numerous results in this direction.
There are several different measures of robustness mentioned in~\cite{Sudakov-survey}, such as random subgraphs, Maker-Breaker games, incompatibility systems and so on.
In this paper the measure of robustness we are interested in is random subgraphs.
For $n\in \mathbb{N}$ and $0\le p\le 1$, we use $G(n,p)$ to denote the Erd\H os--Re\'nyi binomial random graph, which contains $n$ vertices and each pair of vertices is connected by an edge with probability $p$ independent of others.
Given an $n$-vertex graph $G$, we call $G(p):=G\cap G(n,p)$ a \textit{random subgraph} (or \textit{random sparsification}) of $G$. We say that $G(p)$ has a graph property $\mathcal{P}$ \textit{with high probability}, or \textit{w.h.p.} for brevity, if the probability that $G(p)$ has $\mathcal{P}$ tends to 1 as $n$ goes to infinity.
In particular, $G(p)=G(n,p)$ when $G$ is the complete graph $K_n$. The corresponding definitions are made analogously for hypergraphs.

The first result on robustness of graph properties with respect to random subgraphs is a robust version of Dirac's theorem, obtained by Krivelevich, Lee and Sudakov~\cite{MR3180741}.
More precisely, they proved that there exists a constant $C$ such that for $p\ge C\log n /n$ and an $n$-vertex graph $G$ with $\delta(G)\ge n/2$, w.h.p.~$G(p)$ contains a Hamilton cycle.
Recently, Allen, B\"ottcher, Corsten, Davies, Jenssen, Morris, Roberts and Skokan~\cite{Allen}, and independently Pham, Sah, Sawhney and Simkin~\cite{Pham2} proved a robust version of Theorem~\ref{thm:HSz}.

\begin{theorem}[\cite{Allen, Pham2}]\label{thm:robust HSz}
For any integer $r\ge 3$, there exists a constant $C=C(r)$ such that for any $n\in r\mathbb{N}$ and $p\ge Cn^{-2/r}(\log n)^{1/{r \choose 2}}$ the following holds. If $G$ is an $n$-vertex graph with $\delta(G)\ge \left(1-1/r\right)n$, then w.h.p.~$G(p)$ contains a $K_r$-factor.
\end{theorem}

{
The bound on $p$ in Theorem~\ref{thm:robust HSz} is asymptotically tight, as if $p= o(n^{-2/r}(\log n)^{1/{r \choose 2}})$ then w.h.p.~$K_n(p)$ contains vertices not in any copy of $K_r$.
Indeed, the bound on $p$ matches the threshold for $K_r$-factors in random graphs, famously determined by Johansson, Kahn and Vu~\cite{Johansson}.
Therefore, Theorem~\ref{thm:robust HSz} strengthens both Theorem~\ref{thm:HSz} and the Johansson--Kahn--Vu result simultaneously.
}

Inspired by these breakthrough results, Kelly, M\"{u}yesser and Pokrovskiy~\cite{KMP}, and independently Joos, Lang and Sanhueza-Matamala~\cite{Joos-Lang} took a step further and extended Theorem~\ref{thm:robust HSz} to general $F$-factors (and other spanning subgraphs).
To state their result, we need some definitions and notation.
A \textit{$k$-uniform hypergraph} (in short, \textit{$k$-graph}) consists of a vertex set $V$ and an edge set $E \< \binom{V}{k}$, that is, every edge is a $k$-element subset of $V$.
For a $k$-graph $H$ and $d\in [k-1]$, let $\delta_{d}(H)$ be the minimum number of edges of $H$ containing any $d$-element subset of $V(H)$.

\begin{defn}(Dirac threshold)
Let $F$ be a $k$-graph on $s$ vertices.
By $\delta_{F, d}$ we denote the infimum of $\delta$ such that for all $\alpha> 0$ and sufficiently large $n$ the following holds.
Let $H$ be any $k$-graph on $n$ vertices with $n\in s\mathbb{N}$.
If $\delta_{d}(H)\ge (\delta+\alpha)\binom{n-d}{k-d} $, then $H$ contains an $F$-factor.
\end{defn}

For a hypergraph $H$ with at least two vertices, the \textit{$1$-density} of $H$ is $d_{1}(H):=\frac{|E(H)|}{|V(H)|-1}$.
We call $H$ \textit{$1$-balanced }if $d_{1}(H')\le d_{1}(H)$ for every $H'\< H$ and \textit{strictly $1$-balanced} if $d_{1}(H')< d_{1}(H)$ for every $H'\subsetneq H$.

\begin{theorem}[\cite{Joos-Lang}, \cite{KMP}]\label{thm:1-bal}
For every $\alpha>0$ and $k,s\in\N$, there exists $C=C(k,s,\alpha)>0$ such that the following holds for every $d\in[k-1]$. Let $F$ be a strictly $1$-balanced $k$-graph on $s$ vertices, and $H$ be an $n$-vertex $k$-graph such that $\delta_d(H)\ge (\delta_{F,d}+\alpha)\binom{n-d}{k-d}$. If $n\in s\mathbb{N}$ and $p\ge Cn^{-1/d_1(F)}(\log n)^{1/|E(F)|}$, then w.h.p.~$H(p)$ contains an $F$-factor.
\end{theorem}

Again, for strictly $1$-balanced graphs $F$, the bound on $p$ in Theorem~\ref{thm:1-bal} matches the threshold for $F$-factors in random graphs, determined also by Johansson, Kahn and Vu~\cite{Johansson}.

\subsection{Transversal (robust) versions}

Recently, there has been much research on the study of \textit{transversal} versions of Dirac-type results~\cite{MR4125343,MR3628907,MR4394673,rainbow-bandwidth,CHWW1,CHWW2,MR4287703,Ferber,MR4171383,MR4451911,MR4523452,MR4275007,LLL,MR4451150}.
A \textit{graph system} $\mathcal{G}=\{G_1,\ldots,G_m\}$ is a collection of not necessarily distinct graphs on the same vertex set $V$
and an $m$-edge graph $H$ defined on $V$ is \textit{transversal} if $|E(H)\cap E(G_i)|=1$ for each $i\in [m]$. 
The \textit{minimum degree} of a graph system $\mathcal{G}=\{G_1,\ldots,G_m\}$ is defined as the minimum of $\delta(G_i)$ for $i\in [m]$.
The corresponding definitions are made analogously for hypergraphs.
For $F$-factors, Cheng, Han, Wang and Wang~\cite{CHWW1}, and independently Montgomery, M\"uyesser and Pehova~\cite{MR4451150} determined the asymptotic minimum degree condition forcing transversal $K_r$-factors, thus establishing an analogue of Theorem~\ref{thm:HSz} for graph systems, and the latter team indeed obtained essentially optimal results on transversal $F$-factors for arbitrary graph $F$. 

In this paper, we study $F$-factors in the random sparsification of $k$-graph systems, where $F$ is a strictly 1-balanced $k$-graph.
That is, given a $k$-graph system $\mathbf{H} =\{H_{1}, \dots ,H_{m}\}$, we take independent random sparsification of each $k$-graph and consider $\mathbf{H}(p) =\{H_{1}\cap G_1(n,p), \dots ,H_{m}\cap G_m(n,p)\}$.
Although we will also abbreviate each $H_i\cap G_i(n,p)$ as $H_i(p)$, we emphasize here that the sparsifications are done independently (for the case the $k$-graphs in the system are sparsified by the same copy of $G(n,p)$, see Question~\ref{question} in Section 5).
To state our result, let us first define the Dirac threshold in hypergraph systems.
Let $\mathbf{H} =\{H_{1}, H_{2}, \dots ,H_{m}\}$ be a $k$-graph system and $\delta_{d}(\mathbf{H})=\min \{ \delta_{d}(H_{i}) \ \text{for each}\ i\in [m]\}$.

\begin{defn}(Dirac threshold for hypergraph systems)
Let $F$ be a $k$-graph on $s$ vertices with $t$ edges.
By $\delta_{F, d}^{T}$ we denote the infimum of the real number $\delta$ such that for all $\alpha> 0$ and sufficiently large $n$ the following holds.
Let $\mathbf{H} =\{H_{1}, H_{2}, \dots ,H_{tn}\}$ be any $k$-graph system on $sn$ vertices.
If $\delta_{d}(\mathbf{H})\ge (\delta+\alpha)\binom{sn-d}{k-d} $, then $\mathbf{H}$ contains a transversal $F$-factor.
\end{defn}
Now we are ready to state our main result, which is a transversal version of Theorem~\ref{thm:1-bal}.
\begin{theorem}\label{thm: trans 1-bal}
For every $\alpha>0$ and $k,s,t\in\N$, there exists $C=C(k,s,t,\alpha)>0$ such that the following holds for every $d\in[k-1]$. 
Let $F$ be a strictly $1$-balanced $k$-graph on $s$ vertices with $t$ edges, and $\mathbf{H} =\{H_{1}, H_{2}, \dots ,H_{tn}\}$ be any $k$-graph system on $sn$ vertices with $\delta_{d}(\mathbf{H})\ge (\delta^T_{F,d}+\alpha)\binom{sn-d}{k-d} $. If $p\ge Cn^{-1/d_1(F)-1}(\log n)^{1/t}$, then w.h.p.~$\mathbf{H}(p) :=\{H_{1}(p), H_{2}(p), \dots ,H_{tn}(p)\}$ contains a transversal $F$-factor.
\end{theorem}

When $k=2$, by a result of Montgomery, M\"uyesser and Pehova~\cite{MR4451150}, if $F$ has a cut edge or $\delta_{F,1}\ge \frac{1}{2}$, then $\delta_{F,1}^T=\delta_{F,1}$, and thus in this case our result extends Theorem~\ref{thm:1-bal} to graph systems.
On the other hand, for a strictly 1-balanced graph $F$, Theorem~\ref{thm: trans 1-bal} also extends a result of~\cite{MR4451150} on transversal $F$-factors, whose statement is just the $p=1$ case of Theorem~\ref{thm: trans 1-bal}, with $\delta_{F,1}^T$ determined by~\cite{MR4451150}.

For a strictly 1-balanced $k$-graph $F$, Johansson, Kahn and Vu~\cite{Johansson} proved that the threshold for a binomial random $k$-graph $\mathcal{G}^{(k)}(n, p)$ containing an $F$-factor is $\Theta (n^{-1/d_1(F)}(\log n)^{1/|E(F)|})$. 
We note that the lower bound $p\ge Cn^{-1/d_1(F)-1}(\log n)^{1/|E(F)|}$ in Theorem~\ref{thm: trans 1-bal} is best possible up to the constant $C$. Indeed, suppose all graphs in $\mathbf{H}$ are the complete $k$-graph and $p=o(n^{-1/d_1(F)-1}(\log n)^{1/|E(F)|})$. Then using a sprinkling argument, we can regard $\bigcup_{i\in [tn]}H_i(p)$ as the random graph $\mathcal{G}^{(k)}(n,tnp)$. It follows by a well-known result of Rinci\'nski~\cite{Rucinski1992} that not every vertex is covered by a copy of $F$ in $\mathcal{G}^{(k)}(n,tnp)$ as $tnp=o(n^{-1/d_1(F)}(\log n)^{1/|E(F)|})$. 
Last but not least, for general $F$, the situation becomes unclear even in a single binomial random $k$-graph. Let $m_{1}(F)=\max_{F'\< F: |V(F')|>1}d_{1}(F')$. We give an upper bound $O(n^{-1/m_1(F)-1}\log n)$ on the threshold for the existence of transversal $F$-factor.
\begin{theorem}\label{generalF}
For every $\alpha>0$ and $k,s,t\in\N$, there exists $C=C(k,s,t,\alpha)>0$ such that the following holds for every $d\in[k-1]$. 
Let $F$ be a $k$-graph on $s$ vertices with $t$ edges, and $\mathbf{H} =\{H_{1}, H_{2}, \dots ,H_{tn}\}$ be any $k$-graph system on $sn$ vertices with $\delta_{d}(\mathbf{H})\ge (\delta^T_{F,d}+\alpha)\binom{sn-d}{k-d} $. If $p\ge Cn^{-1/m_1(F)-1}\log n$, then w.h.p.~$\mathbf{H}(p) :=\{H_{1}(p), H_{2}(p), \dots ,H_{tn}(p)\}$ contains a transversal $F$-factor.
\end{theorem}


\subsection{Related works}
In recent years, numerous results have emerged on the robustness of Dirac-type theorems under random sparsifications, with the robust thresholds for various spanning structures being determined.
Kang,  Kelly, K\"uhn, Osthus and Pfenninger~\cite{Kang-Kelly} and Pham, Sah, Sawhney and Simkin~\cite{Pham2} independently  determined the robust threshold for 
perfect matchings to be $\Omega(\log n/n^{k-1})$.
In~\cite{Pham2}, they also proved that for any $n$-vertex graph $G$ with $\de(G)\ge (1/2+\alpha)n$, w.h.p.~$G(C\log n/n)$ contains a given bounded-degree spanning tree.
Later, Bastide, Legrand-Duchesne and M\"uyesser~\cite{Bastide} provided an alternative proof of this result. 
As mentioned earlier, the authors in~\cite{Joos-Lang} and~\cite{KMP} independently determined the robust threshold for 
$F$-factors, where $F$ is strictly 1-balanced. Furthermore, they also independently considered robust thresholds for Hamiltonicity, including power of Hamilton cycles, loose and tight Hamilton cycles.

For the transversal version of Dirac-type results, Joos and Kim~\cite{MR4171383} proved a transversal analogue of Dirac’s theorem in 2020, confirming a conjecture of Aharoni (see~\cite{MR4125343}) and improving upon an asymptotically tight result of Cheng, Wang and Zhao~\cite{MR4287703}. 
More precisely, they~\cite{MR4171383} showed that any graph system $\mathcal{G}=\{G_1,\ldots, G_n\}$ on $n$ vertices with $\de(\mathcal{G})\ge n/2$ contains a transversal Hamilton cycle. By taking $G_1=\cdots=G_n$, the result recovers Dirac’s theorem and thus generalizes it.
Since then, there has been substantial progress in the study of the transversal version of Dirac-type results. We refer interested readers to a comprehensive survey by Sun, Wang and Wei~\cite{SWW}.

Recently, Ferber, Han and Mao~\cite{Ferber} derived a transversal robust version of Dirac's theorem by 
showing that for any graph system $\mathcal{G}=\{G_1,\ldots, G_n\}$ on $n$ vertices with $\de(\mathcal{G})\ge (1/2+\alpha)n$, if $p=\Omega(\log n/n)$, then w.h.p.~$\{G_1(p),\ldots,G_n(p)\}$ contains a transversal Hamilton cycle.
In fact, they proved a transversal version of resilience result on Hamiltonicity which implies the above result. Very recently, Anastos and Chakraborti~\cite{Anastos-Chakraborti} improved the above result by removing the error term in the minimum degree condition with spread techniques. 
Similar result for transversal perfect matchings is also obtained in~\cite{Ferber} by using resilience results on perfect matchings in random bipartite graphs due to Sudakov and Vu~\cite{Sudakov2008}.


\subsection{Basic notation}
For positive integers $a,b$ with $a<b$, let $[a]=\{1,2,\ldots,a\}$ and $[a,b]=\{a,a+1,\ldots,b\}$.
For constants $x,y,z$, $x=y\pm z$ means that $y-z\le x\le y+z$, and $x\ll y$ means that for any $y> 0$ there
exists $x_0> 0$ such that for any $x< x_0$ the subsequent statement holds.


For hypergraphs $F$ and $H$, the \textit{$F$-complex} of $H$, denoted by $H_F$, is the $|V(F)|$-uniform multi-hypergraph with vertex set $V(H)$ in which every copy of $F$ in $H$ corresponds to a distinct hyperedge of $H_F$ on the same set of vertices. Note that $H$ has an $F$-factor if and only if $H_F$ contains a perfect matching. 
Let $\mathcal{G}^{(k)}(n, p)$ be the binomial random $k$-graph on $n$ vertices where every edge of  the complete $k$-graph is included independently with probability $p$.
Let $\mathcal{G}_F(n, p)$ be the $F$-complex of $\mathcal{G}^{(k)}(n, p)$.

For a (hyper)graph $G$, we use $e(G)$ to denote the number of edges in $G$. Given a hypergraph system $\mathbf{H}=\{H_{1}, H_{2},\dots,H_{m}\}$, let $V(\mathbf{H})$ denote the vertex set of $H_{i}$ for each $i\in [m]$.
For each $i\in [m]$, we can assume that the edges of $H_{i}$ is colored with respect to color $i$.
Let $\mathcal{C}(\mathbf{H})=[m]$ denote the color set of $\mathbf{H}$.

\section{Spreadness}


The \textit{spread method} (also known as the \textit{fragmentation method}) has recently brought major breakthroughs in random graph theory, featured in the celebrated resolution of the \textit{fractional expectation threshold versus threshold} conjecture of Talagrand~\cite{Talagrand}, by Frankston, Kahn, Narayanan and Park~\cite{Frankston} (see Theorem~\ref{lem:Frankston} below). 
It is worth mentioning that Park and Pham~\cite{Park-Pham} proved the Kahn--Kalai conjecture~\cite{Kahn-Kalai}, which implies Talagrand's conjecture.
By constructing probability measures with good spread, the spread method has powerful applications to various difficult problems, such as robustness and enumeration results (see e.g. \cite{Anastos-Chakraborti, Bastide, Joos-Lang, Kang-Kelly, KMP, Pham2}).

Given (hyper)graphs $H$ and $G$, we will consider a hypergraph $(V, \mathcal{G})$ where $V:=E(H)$ and $\mathcal{G}$ is the family of all subgraphs of $H$ which are isomorphic to $G$.
Suppose that $e(G)=r$. Then $(V, \mathcal{G})=(E(H), \mathcal{G})$ is an $r$-graph.
\begin{defn}[Spread]
Let $q\in [0,1]$ and $r\in \mathbb{N}$. Let $(E(H),\mathcal{G})$ be an $r$-graph as above, and let $\mu$ be a probability distribution on $\mathcal{G}$. 
We say  $\mu$ is \textit{$q$-spread} if
\begin{align*}
\mu (\{A\in \mathcal{G}: A\supseteq S\})\le q^{|S|} \ \text{for all}\ S\< E(H).
\end{align*}
\end{defn}

For example, when $G$ is a perfect matching on $n$ vertices and $H=K_n$, we consider the uniform measure $\mu$ supported on all perfect matchings in $H$. Using Stirling's formula one can get that $\mu$ is actually $(e/n)$-spread where $e$ is the base of natural logarithm. 
A recent result (among others) of Pham, Sah, Sawhney, and Simkin~\cite{Pham2} studies the spread distribution on perfect matchings in $\eps$-super regular graphs as well as in dense graphs as follows. 
A bipartite graph $H=(V_1,V_2)$ is \textit{balanced} if $|V_1|=|V_2|$.
\begin{theorem}[\cite{Pham2}]\label{thm:bipartite spread}
There exists $C>0$ with the following property.
If $H$ is a balanced bipartite graph on $2n$ vertices with $\delta(H)\ge 3n/4$, then there exists a $(C/n)$-spread distribution on perfect matchings in $H$.
\end{theorem}

Another crucial notion, vertex-spread, was introduced in~\cite{Pham2} and has been successful in this line of research.
A \textit{hypergraph embedding} $\varphi: G\hookrightarrow H$ of a hypergraph $G$ into a hypergraph $H$ is indeed an injection $\varphi: V(G)\rightarrow V(H)$ that maps edges of $G$ to edges of $H$.
Thus, there is an embedding of $G$ into $H$ if and only if $H$ contains a subgraph isomorphic to $G$.

\begin{defn}[Vertex-spread]
Let $G$ and $H$ be two finite hypergraphs, and let $\mu$ be a probability distribution on all embeddings $\varphi : G \hookrightarrow H$.
For $q\in [0,1]$, we say that $\mu$ is a \textit{$q$-vertex-spread} measure if for every two sequences of distinct vertices $x_{1}, \dots , x_{s}\in V(G)$ and $y_{1}, \dots , y_{s}\in V(H)$,
\begin{align*}
\mu(\{\varphi: \varphi(x_{i})=y_{i} \ \text{for all}\ i\in [s]\})\le q^{s}.
\end{align*}
\end{defn}
Unlike the spread distribution which concerns the probability that a prescribed set of edges (in $H$) are contained in the random sample of hyperedges, vertex-spreadness measures the probability that a set of vertices in $G$ are mapped into the same number of prescribed vertices in $H$ in a random embedding $G\hookrightarrow H$. The following result shows that the vertex-spread distribution can be transformed into a spread distribution. Recall that $m_{1}(G)=\max_{G'\< G: |V(G')|>1}d_{1}(G')$, where $d_{1}(G'):=\frac{|E(G')|}{|V(G')|-1}$.

\begin{prop}[\cite{KMP}]\label{prop: vtx-spread and spread}
For every $C, k, \De >0$, there exists $C'=C'(C, k, \De)>0$ such that the following holds for all sufficiently large $n$. Let $H$ and $G$ be $n$-vertex $k$-graphs. If there is a $(C/n)$-vertex-spread distribution on embeddings $G\hookrightarrow H$ and $\De(G)\le \De$, then there is a $(C'/n^{1/m_{1}(G)})$-spread distribution on subgraphs of $H$ which are isomorphic to $G$. \end{prop}

As the first attempt, we extend Theorem~\ref{thm:bipartite spread} to a vertex-spread version as follows (see Theorem~\ref{thm:bipartite perfect matching} for a stronger statement). 

\begin{theorem}\label{thm:bipartite perfect}
For every $\varepsilon > 0$ there exists $C=C(\varepsilon )>0$ such that the following holds for every sufficiently large $n$.
Let $G=(V_{1}, V_{2})$ be a balanced bipartite graph such that $|V_{1}|=|V_{2}|=n$.
If $\delta(G)\ge (\tfrac{1}{2}+\varepsilon )n$, 
then there is a $(C/n)$-vertex-spread distribution on embeddings of perfect matchings in $G$.  
\end{theorem}

Combined with Proposition~\ref{prop: vtx-spread and spread}, our result immediately gives a $O(1/n)$-spread distribution over perfect matchings in a balanced bipartite graph on $2n$ vertices with minimum degree at least $(\tfrac{1}{2}+o(1))n$.



\subsection{Vertex-spread of perfect matchings in $k$-partite $k$-graphs}

It is well-known that the decision problem asking whether there exists a perfect matching in $3$-partite $3$-graphs is one of the first $21$ NP-Complete problems. Over the last few decades,
the question of determining the minimum degree threshold forcing a perfect matching in $k$-partite $k$-graphs has received considerable attention from extremal graph theory. A multitude of cornerstone results have been presented on determining the various minimum degree thresholds for perfect matchings in a series of works~\cite{Aharoni,  HPS, DKDO, Pikhurko, Rodl}.
To formulate this, for a $k$-partite $k$-graph $H$ with vertex parts $V_{1}\cup \dots \cup V_{k}= V(H)$ and an index set $L\subseteq[k]$, let $\delta_{L}(H)$ be the minimum of $|N(X)|$ over all sets $X\< V(H)$ such that $|X|=|L|$ and $|X\cap V_{i}|=1$ for each $i\in L$.
In particular, for a bipartite graph $G=(A, B)$, we denote by $\delta_{\{1\}}(G)$ the smallest degree of a vertex from $A$ and by $\delta_{\{2\}}(G)$ the smallest degree of a vertex from $B$. 

Aharoni, Georgakopoulos and Spr{\"u}ssel~\cite{Aharoni} proved that the minimum degree condition $\delta_{k-1}(H)\ge n/2$ can guarantee the existence of a perfect matching in a balanced $k$-partite $k$-graph $H$. 
Pikhurko~\cite{Pikhurko} considered the same problem under a `cooperative' minimum degree condition.

\begin{theorem}[\cite{Pikhurko}]\label{thm:multipartite perfect matching}
Let $k\ge 2$, $\ell\in [k-1]$, and $L\in \binom{[k]}{\ell}$ be fixed. Let $n$ be sufficiently large and
\begin{align*}
\lambda=\sqrt{257kn\log n}.
\end{align*}
Let $H$ be a $k$-partite $k$-graph with parts $V_{1}\cup \dots \cup V_{k}=V$ such that $|V_{i}|=n$ for each $i\in [k]$.
If
\begin{align*}
\delta_{L}(H)n^{\ell}+\delta_{[k]\setminus L}(H)n^{k-\ell}\ge n^{k}+k\lambda n^{k-1},
\end{align*}
then $H$ admits a perfect matching.
\end{theorem}
Following Pikhurko's strategy, we can actually derive a $O(1/n)$-vertex-spread distribution on embeddings of perfect matchings into a balanced $k$-partite $k$-graph. 

\begin{theorem}\label{thm:multipartite spread}
For every $\varepsilon > 0$ and $k \in \mathbb{N}$ there exists $C=C(\varepsilon, k)>0$ such that the following holds for every $\ell \in [k-1]$ and sufficiently large $n$.
Let $L\in \binom{[k]}{\ell}$ and
$H$ be a $k$-partite $k$-graph with parts $V_{1}\cup \dots \cup V_{k}=V(H)$ such that $|V_{i}|=n$ for each $i\in [k]$.
If
\begin{align*}
\delta_{L}(H)n^{\ell}+\delta_{[k]\setminus L}(H)n^{k-\ell}\ge (1+\varepsilon )n^{k},
\end{align*}
then there is a $(C/n)$-vertex-spread distribution on embeddings of perfect matchings into $H$.
\end{theorem}

To achieve this, we actually first prove the basic bipartite case in the following, that is, to construct a desired vertex-spread distribution on embeddings of perfect matchings in bipartite graphs. 

\begin{theorem}\label{thm:bipartite perfect matching}
For every $\varepsilon > 0$ there exists $C=C(\eps)>0$ such that the following holds for every sufficiently large $n$.
Let $G=(V_{1}, V_{2})$ be a balanced bipartite graph such that $|V_{1}|=|V_{2}|=n$.
If
\begin{align*}
\delta_{\{1\}}(G)+\delta_{\{2\}}(G)\ge (1+\varepsilon )n,
\end{align*}
then there is a $(C/n)$-vertex-spread distribution on embeddings 
of perfect matchings in $G$.  
\end{theorem}
Therefore, combining this with Pikhurko's strategy, we can easily prove Theorem~\ref{thm:multipartite spread}, whose proof is deferred to Appendix~\ref{app1}.

\subsection{Vertex-spread of general $F$-factors}
For more general $F$-factors (rather than perfect matchings), a recent significant insight by Kelly, M\"uyesser and Pokrovskiy~\cite{KMP} tells that one can always derive a $O(1/n)$-vertex-spread distribution on all embeddings of $F$-factors whenever the minimum degree of the host graph slightly exceeds the threshold for the containment of an $F$-factor.

\begin{theorem}[Theorem 1.15 in \cite{KMP}]
For every $\alpha>0$ and $k, r\in \mathbb{N}$ there exists $C=C(\alpha, k, r)>0$ such that the following holds for all $d\in [k-1]$ and all sufficiently large $n$ for which $r$ divides $n$. Let $F$ be an $r$-vertex $k$-graph. If $H$ is an $n$-vertex $k$-graph such that $\de_d(H)\ge (\de_{F,d}+\alpha)\binom{n-d}{k-d}$, then there is a $(C/n)$-vertex-spread distribution on the embeddings of $F$-factors in $H$.
\end{theorem}

We shall generalize the vertex-spread result as above to the transversal setting.
To state our result, we need to make some preparations. 
Let $F$ be a $k$-graph on $s$ vertices with $t$ edges and $\mathbf{H}=\{H_{1}, H_{2}, \dots , H_{tn}\}$ be a $k$-graph system on $sn$ vertices.
We first construct a $(k+1)$-graph $\mathcal{H}$ based on $\mathbf{H}$. 
The vertex set of $\mathcal{H}$ is $V(\mathbf{H})\cup \mathcal{C}(\mathbf{H})$.
Every edge of $\mathcal{H}$ consists of $k$ elements in $V(\mathbf{H})$ and one element in $\mathcal{C}(\mathbf{H})$. For any $\{v_1,\ldots,v_k\}\< V(\mathbf{H})$ and $i\in \mathcal{C}(\mathbf{H})$, $\{v_1,\ldots,v_k,i\}$ is an edge of $\mathcal{H}$ if and only if $\{v_1,\ldots,v_k\}$ is an edge of $H_i$.
In this case, we define the \textit{expansion} of $F$, denoted by $F^{*}$, which is obtained from $F$ by adding an extra element to each edge of $F$ such that different edges of $F$ receive different elements. Hence $F^{*}$ is a $(k+1)$-graph on $(s+t)$ vertices with $t$ edges. Here we may abuse the notation and call these extra elements the \textit{colors} in $F^*$.  Note that a transversal $F$-factor in $\mathbf{H}$ is equivalent to an $F^{*}$-factor in $\mathcal{H}$.

\begin{theorem}\label{thm:hypergraph system spread}
For every $\alpha > 0$ and $k, s, t \in \mathbb{N}$ there exists $C=C(\alpha, k, s, t)>0$ such that the following holds for every $d\in [k-1]$ and every sufficiently large $n$.
Let $F$ be a $k$-graph on $s$ vertices with $t$ edges, and let $\mathbf{H}=\{H_{1}, H_{2}, \dots , H_{tn}\}$ be a $k$-graph system on $sn$ vertices.
If 
\begin{align*}
\delta_{d}(H_{i})\ge (\delta_{F,d}^{T}+\alpha)\binom{sn-d}{k-d} \ \text{for each} \ i\in [tn],
\end{align*}
then there is a $(C/n)$-vertex-spread distribution on the embeddings of $F^{*}$-factors into $\mathcal{H}$, where $F^{*}$ and $\mathcal{H}$ are defined as above.
\end{theorem}
An immediate corollary of Theorem~\ref{thm:hypergraph system spread} is that the hypergraph system $\mathbf{H}$ contains  $e^{(s+t-1)n\log n-O(n)}$ distinct transversal $F$-factors. 
Indeed, if there is a $(C/n)$-vertex-spread distribution $\mu$ on embeddings of $F^{*}$-factors in $\mathcal{H}$, then we have
\begin{align*}
1=\sum_{\varphi: G \hookrightarrow \mathcal{H}}\mu(\{\varphi\})\le |\{\varphi: G \hookrightarrow \mathcal{H}\}|\cdot\left(\frac{C}{n}\right)^{sn+tn}=|\{\varphi: G \hookrightarrow \mathcal{H}\}|\cdot e^{(s+t)n\log C-(s+t)n\log n},
\end{align*}
where $\varphi$ is taken over all embeddings of $F^{*}$-factors in $\mathcal{H}$.
Note that the automorphism number of an $F^{*}$-factor is $n!\cdot(aut(F^{*}))^n=e^{n\log n-O(n)}$. Hence there are $e^{(s+t-1)n\log n-O(n)}$ distinct $F^{*}$-factors in $\mathcal{H}$ and so there are $e^{(s+t-1)n\log n-O(n)}$ distinct transversal $F$-factors in $\mathbf{H}$.

\subsection{Proof of the main theorem via the spread method}
In this section, we shall show how to deduce our main result, Theorem~\ref{thm: trans 1-bal}, from Theorem~\ref{thm:hypergraph system spread}.
For this, we also need the following two well-known results. The first one is the aforementioned breakthrough of Frankston, Kahn, Naranyanan, and Park~\cite{Frankston} on the fractional expectation threshold conjecture of Talagrand~\cite{Talagrand}.
The following theorem can be deduced from Theorem 1.6 of \cite{Frankston}.

\begin{theorem}[\cite{Frankston}]\label{lem:Frankston}
There exists $C>0$ such that the following holds. Let $X$ be a finite set, and $\mathcal{F}\< 2^X$. If there is a $q$-spread probability measure on $\mathcal{F}$, then w.h.p.~$X(\min\{Cq\log |X|,1\})$ has a subset which is an element of $\mathcal{F}$.
\end{theorem} 

The second result we will use is a coupling result of Riordan~\cite{Riordan}.
\begin{theorem}[\cite{Riordan}]\label{thm: coupling-hypergraph}
For every $k, r \in \mathbb{N}$, there exists $c \in (0,1]$ such that the following holds. If $F$ is a strictly 1-balanced $k$-graph on $r$ vertices and $p = p(n) \leq \log^2(n) / n^{1 / d_1(F)}$, then for some $\pi = \pi(n) \sim cp^{e(F)}$, we may couple $G \sim \mathcal{G}^{(k)}(n, p)$ with $G_F \sim \mathcal{G}_F(n, \pi)$ such that, w.h.p.~for every $F$-edge present in $G_F$, the corresponding copy of $F$ is present in $G$.
\end{theorem}

Now we are ready to prove Theorem~\ref{thm: trans 1-bal} and Theorem~\ref{generalF}.

\begin{proof}[Proof of Theorem~\ref{thm: trans 1-bal}]
We reformulate our original problem into finding an $F^{*}$-factor in $\mathcal{H}(p)$, where $F^{*}$ and $\mathcal{H}$ are defined just above Theorem~\ref{thm:hypergraph system spread}.
Let $\mathcal{H}_{F^{*}}$ be the $F^{*}$-complex of $\mathcal{H}$ and so $\mathcal{H}_{F^{*}}$ is $(s+t)$-uniform. Note that a perfect matching in $\mathcal{H}_{F^{*}}$ is equivalent to an $F^{*}$-factor in $\mathcal{H}$, and is also equivalent to a transversal $F$-factor in $\mathbf{H}$.

\begin{claim}\label{claim: strictly 1-balanced}
    $F^{*}$ is strictly 1-balanced\footnote{For general graph $F$, one can similarly show that $m_1(F^*)=\frac{m_1(F)}{1+m_1(F)}$.}. 
\end{claim} 
\begin{proof}[Proof of Claim~\ref{claim: strictly 1-balanced}]
Let $F'$ be any proper subgraph of $F^{*}$. Assume that $F'$ has $x$ vertices from $[tn]$ and $y$ vertices from $V(F)$ which induce $y'$ edges of $F$. Since $F$ is strictly 1-balanced, we have $y'/(y-1)< t/(s-1)$ and thus $y-1>y'(s-1)/t$. Note that $d_1(F^{*})=t/(s+t-1)$ and then
\[d_1(F')\le \tfrac{\min\{x,y'\}}{x+y-1}<\tfrac{\min\{x,y'\}}{x+\tfrac{s-1}{t}y'}\le \tfrac{t}{s+t-1}=d_1(F^{*}),\] 
where the last inequality follows as $\min\{x,y'\}\le \frac{tx+(s-1)y'}{t+s-1}$. Hence $F^{*}$ is strictly 1-balanced.
\end{proof}

By Theorem~\ref{thm:hypergraph system spread}, there is a $(C_{\ref{thm:hypergraph system spread}}/n)$-vertex-spread distribution on the embeddings of $F^{*}$-factors into $\mathcal{H}$. This immediately gives a $(K/n)$-vertex-spread distribution on embeddings of perfect matchings into $\mathcal{H}_{F^{*}}$ for some $K=K(C_{\ref{thm:hypergraph system spread}},F)$. 
Then by Proposition~\ref{prop: vtx-spread and spread},  
there is a $\left( C_{\ref{prop: vtx-spread and spread}}n^{-(s+t-1)} \right)$-spread measure on perfect matchings of $\mathcal{H}_{F^{*}}$.
By applying Theorem~\ref{lem:Frankston} with $X$ being the edge set of $\mathcal{H}_{F^{*}}$ and $\mathcal{F}$ being the family of perfect matchings in $\mathcal{H}_{F^{*}}$, we obtain that w.h.p.~$\mathcal{H}_{F^{*}}(\pi)$ has a perfect matching, where $\pi=C_{\ref{lem:Frankston}}n^{-(s+t-1)}\log n$.
Then Theorem~\ref{thm: coupling-hypergraph}, combined with Claim~\ref{claim: strictly 1-balanced}, implies that w.h.p.~$\mathcal{H}(p)$ contains an $F^{*}$-factor, where $p\ge p_0=C(n^{-(s+t-1)}\log n)^{1/t}=Cn^{-1/d_1(F)-1}(\log n)^{1/t}$ for some constant $C$.
That is, w.h.p.~$\mathbf{H}(p) :=\{H_{1}(p), H_{2}(p), \dots ,H_{tn}(p)\}$ contains a transversal $F$-factor.
\end{proof}

The proof of Theorem~\ref{generalF}
closely follows that of Theorem~\ref{thm: trans 1-bal} and the only difference is discussed as follows. From Theorem~\ref{thm:hypergraph system spread}, there is a $(C_{\ref{thm:hypergraph system spread}}/n)$-vertex-spread distribution on the embeddings of $F^{*}$-factors into $\mathcal{H}$. In this case, Proposition~\ref{prop: vtx-spread and spread} immediately gives a $(C_{\ref{prop: vtx-spread and spread}}/n^{1/m_1(F^*)})$-spread measure on $F^{*}$-factors in $\mathcal{H}$. By applying Theorem~\ref{lem:Frankston} with $X$ being the edge set of $\mathcal{H}$ and $\mathcal{F}$ being the family of $F^{*}$-factors in $\mathcal{H}$, we obtain that for some $C>0$ and $p\ge Cn^{-1/m_1(F^*)}\log n=Cn^{-1/m_1(F)-1}\log n$, w.h.p.~$\mathcal{H}(p)$ contains a $F^{*}$-factor as $m_1(F^*)=\frac{m_1(F)}{1+m_1(F)}$.
\section{Proof ideas: Random clustering method}

In this section, we shall present some essential tools used in the proofs of Theorems~\ref{thm:bipartite perfect matching} and~\ref{thm:hypergraph system spread}. 
We employ the random clustering strategy recently introduced by Kelly, M\"uyesser and Pokrovskiy~\cite{KMP}. 

\subsection{Proof outlines}
We first review the idea of Kelly, M\"uyesser and Pokrovskiy~\cite{KMP}.
The strategy starts with randomly partitioning the vertex set of the host graph $G$ into clusters of (large) constant sizes. Note that most of the clusters will inherit the minimum degree condition of $G$, but there might be some \textit{exceptional} clusters. 
To overcome this, one may collect the vertices from the exceptional clusters (only a few) and redistribute them to other clusters (which we call \textit{good}), so that the minimum degree conditions can be retained, by building an auxiliary bipartite graph and then finding a spread perfect matching. 
Technically, the redistribution process is based on a $O(1/n)$-spread distribution on perfect matchings as given in Theorem~\ref{thm:bipartite spread}. A similar strategy also appears in a recent subsequent work~\cite{Bastide} on embedding bounded-degree spanning trees. 

For our problem of transversal $F$-factors in Theorem~\ref{thm:hypergraph system spread}, we shall try to randomly split both the vertex set $V(\mathbf{H})$ and the color set  $\mathcal{C}(\mathbf{H})$ into clusters of constant sizes, and consider \textit{sub-systems}, which consist of, e.g., a set of $sC$ vertices and $tC$ colors for some large integer $C>0$. 
Then to proceed the (spread) matching process, 
we shall argue via an auxiliary hypergraph of higher uniformity, as we need to redistribute a small number of both vertices and colors to good clusters.
That is, after we collect the vertices and colors from the exceptional clusters, we may try to redistribute $s$ vertices and $t$ colors to good a cluster, and in this way we shall define and use an auxiliary $(s+t+1)$-partite $(s+t+1)$-graph to guide the matching process.
Then the spread perfect matching can be built by Theorem~\ref{thm:multipartite spread}.

As mentioned in Section 2, we shall prove Theorem~\ref{thm:multipartite spread} following the proof idea of Pikhurko and reduce it to Theorem~\ref{thm:bipartite perfect matching}.
To prove Theorem~\ref{thm:bipartite perfect matching}, consider a balanced bipartite graph $G=(V_1,V_2)$, and we need to embed a perfect matching in $G$ with a $(C/n)$-vertex-spread distribution. 
We also use the random clustering method and first split $V(G)$ into random clusters of (large) constant sizes. 
Note that most of the clusters inherit the minimum degree condition. 
Then we randomly redistribute the vertices of bad clusters into good ones so that the minimum degree condition is maintained. 
Note that as $G$ is bipartite, this is similar to the graph system case -- each good cluster must receive an equal number of vertices from $V_1$ and $V_2$.
To do this, we distribute the vertices in two rounds. 
By building two respective auxiliary bipartite graphs, we ensure that each good cluster receives exactly one vertex from $V_1$ and one vertex from $V_2$. Combining the respective spread distributions of perfect matchings (from Theorem~\ref{thm:bipartite spread}) in the two auxiliary bipartite graphs, we obtain a desired vertex-spread distribution on embeddings of perfect matchings as in Theorem~\ref{thm:bipartite perfect matching}.\medskip

\subsection{Two lemmas}
We now state our two random cluster lemmas, Lemmas~\ref{lem:clustering} and~\ref{lem:bipartite clustering}.
Here is a proof outline.
We will use Lemma~\ref{lem:bipartite clustering} to prove Theorem~\ref{thm:bipartite perfect matching} and as mentioned earlier, Theorem~\ref{thm:bipartite perfect matching} implies Theorem~\ref{thm:multipartite spread}. 
Moreover, Theorem~\ref{thm:multipartite spread} will be used in the proof of Lemma~\ref{lem:clustering}, which in turn is used to prove Theorem~\ref{thm:hypergraph system spread}.

\begin{lemma}\label{lem:clustering}
Let $1/n\ll 1/C' \ll 1/C \ll \alpha,  1/k, 1/s, 1/t$, where $n, k, s, t, C \in \mathbb{N}$, and let $\delta \in [0,1]$.
If $\mathbf{H} =\{H_{1}, H_{2}, \dots ,H_{tn}\}$ is a $k$-graph system with $\delta_{d}(\mathbf{H})\ge (\delta+\alpha)\binom{sn-d}{k-d}$,
then there exists a probability distribution over all partitions $\mathcal{U}:=\{U_{1}, U_{2}, \dots, U_{m}\}$ of $V(\mathbf{H})$ and $\mathcal{W}:=\{W_{1}, W_{2}, \dots, W_{m}\}$ of $\mathcal{C}(\mathbf{H})$ satisfying the following properties:

\begin{enumerate}[label = (\arabic{enumi})]
\rm \item \label{size} $|U_{i}|=sC$, $|W_{i}|=tC$ for each $2\le i \le m$ and $|U_{1}|$ equals $sC(C-1)$ plus the remainder when $sn$ is divided by $sC(C-1)$, $|W_{1}|$ equals $tC(C-1)$ plus the remainder when $tn$ is divided by $tC(C-1)$;
\rm \item \label{degree} for every $i\in [m]$, writing $\mathbf{H}_{i}=\{H_{j}[U_{i}]\mid j\in W_{i}\}$, we have
$\delta_{d}(\mathbf{H}_{i})\ge (\delta+\alpha/2)\binom{|U_{i}|-d}{k-d}$;

\rm \item \label{location} for all distinct vertices $y_{1}, y_{2},\dots , y_{r}\in V(\mathbf{H})$ and colors $c_{1}, c_{2},\dots , c_{\ell }\in \mathcal{C}(\mathbf{H})$, and functions $f: [r]\to [m]$, $g: [\ell]\to [m]$, we have 
\[\mathbb{P}\big{[}y_{i}\in U_{f(i)}~\text{for all}~i\in [r] ~\text{and}~ c_{j}\in W_{g(j)}~\text{for all}~j\in [\ell]\big{]} \le \left(\frac{C'}{n}\right)^{r+\ell}.\]
\end{enumerate}
\end{lemma}


Before stating Lemma~\ref{lem:bipartite clustering}, we give a proof of Theorem~\ref{thm:hypergraph system spread} using Lemma~\ref{lem:clustering}.

\begin{proof}[Proof of Theorem~\ref{thm:hypergraph system spread} from Lemma~\ref{lem:clustering}] 
Let $\alpha>0$, $k \in \mathbb{N}$ and $d\in [k-1]$. Fix $F$ to be a $k$-graph on $s$ vertices with $t$ edges and choose $1/n \ll 1/C' \ll 1/C \ll \alpha,  1/k, 1/s, 1/t$ as in Lemma~\ref{lem:clustering}. 
Given a $k$-graph system $\mathbf{H}=\{H_{1}, H_{2}, \dots , H_{tn}\}$ on $sn$ vertices with $\delta_{d}(\mathbf{H})\ge (\delta_{F,d}^{T}+\alpha)\binom{sn-d}{k-d}$.
Now we embed a transversal $F$-factor into $\mathbf{H}$ in the following two steps.
\begin{itemize}
\item [(i)] Randomly sample partitions $\mathcal{U}:=\{U_{1}, U_{2}, \dots, U_{m}\}$ and $\mathcal{W}:=\{W_{1}, W_{2}, \dots, W_{m}\}$ as required in Lemma~\ref{lem:clustering} with $\delta_{F,d}^{T}$ playing the role of $\delta$.

\item [(ii)] By combining property~\ref{degree} of Lemma~\ref{lem:clustering} and the definition of $\de_{F,d}^T$, we can find a transversal $F$-factor in each sub-system $\mathbf{H}_{i}=\{H_{j}[U_{i}]: j\in W_{i}\}$. For each $i\in [m]$, we choose an embedding of transversal $F$-factor in $\mathbf{H}_{i}$ uniformly at random (among all embeddings), which altogether form an embedding $\phi$ of transversal $F$-factor into $\mathbf{H}$.
\end{itemize}
Recall that $F^*$ is the expansion of $F$ and thus part (ii) gives an $F^*$-factor in the corresponding $\mathcal{H}$. Thus $\phi$ is also an embedding of $F^*$-factor into $\mathcal{H}$. 
Let $G$ be the graph consisting of $n$ vertex-disjoint copies of $F^*$ and use $C(G)$ to denote the set of colors from every copy of $F^*$. Without loss of generality, we may assume $C(G)=[tn]=Q_{1}\cup \dots \cup Q_{m}$ where $Q_{1}=\{1,2,\cdots,|W_1|\}$, $Q_{i}=\{|Q_{i-1}|+1, \cdots, |Q_{i-1}|+|W_i|\}$ and $V(G)\setminus C(G)=[sn]=P_{1}\cup \dots \cup P_{m}$ where $P_{1}=\{1,2,\cdots,|U_1|\}$, $P_{i}=\{|P_{i-1}|+1, \cdots, |P_{i-1}|+|U_i|\}$ for $2\le i\le m$.  Additionally, we assume $G[P_i\cup Q_i]$ induces an $F^*$-factor for each $i\in [m]$.

It remains to show that the random embedding $\phi$ admits a $(C'/n)$-vertex-spread distribution, 
that is, for all distinct sequences \[x_{1}, x_{2},\dots , x_{r}\in V(G)\setminus C(G), ~n_{1}, n_{2},\dots , n_{\ell}\in C(G)\]
\[\text{and}~~ y_{1}, y_{2},\dots , y_{r}\in V(\mathbf{H}),~ c_{1}, c_{2},\dots , c_{\ell }\in \mathcal{C}(\mathbf{H}),\] we have 
$\mathbb{P}[\phi(x_{i})=y_{i} \ \text{and}\ \phi(n_{j})=c_{j}\ \text{for all}\ i\in [r] ,   j\in [\ell]] \le (\frac{C'}{n})^{r+\ell}$.
Define $f: [r]\rightarrow[m]$ and $g: [\ell]\rightarrow[m]$ such that $f(i)$ is the unique index with $x_{i}\in P_{f(i)}$, and $g(j)$ is the unique index with $n_{j}\in Q_{g(j)}$, respectively.
By the definition of $\phi$, we know that
\begin{align*}
&\mathbb{P}\big{[}\phi(x_{i})=y_{i}\ \text{and}\ \phi(n_{j})=c_{j} \ \text{for all}\ i\in [r] , j\in [\ell]\big{]}\\
\le 
&\mathbb{P} \big{[}y_{i}\in U_{f(i)}~\text{and}~c_{j}\in W_{g(j)}~\text{for each}~i\in [r]~\text{and}~ j\in [\ell]\big{]}
\le\left(\frac{C'}{n}\right)^{r+\ell},
\end{align*}
where the last inequality holds by~\ref{location} of Lemma~\ref{lem:clustering}.
\end{proof}

Now we give the statement of the second lemma on random clustering.

\begin{lemma}\label{lem:bipartite clustering}
Let $1/n \ll 1/C' \ll 1/C \ll \varepsilon $, where $n\in \mathbb{N}$.
If $G=(V_{1}, V_{2})$ is a bipartite graph with $|V_{1}|=|V_{2}|=n$ such that $\delta_{\{1\}}(G)+\delta_{\{2\}}(G)\ge (1+\varepsilon )n$,
then there exists a probability distribution over all partitions 
$\mathcal{U}:=\{U_{1}, U_{2}, \dots, U_{m}\}$ of $V(G)$ satisfying the following properties:

\begin{enumerate}[label = (\arabic{enumi})]
\rm \item\label{size2} $|U_i\cap V_1|=|U_i\cap V_2|$ for each $i\in [m]$, $|U_{i}|=2C$ for each $2\le i\le m$ and $|U_{1}|$ equals $2(C-1)C$ plus the remainder when $2n$ is divided by $2(C-1)C$;
\rm \item\label{degree2} $\delta_{\{1\}}(G[U_{i}])+\delta_{\{2\}}(G[U_{i}])\ge (1+\eps/3)\frac{|U_{i}|}{2}$ for each $i\in [m]$;
\rm \item\label{location2}  for every set of distinct vertices $y_{1}, y_{2},\dots , y_{s}\in V(G)$ and every function $f: [s]\to [m]$, we have 
$$\mathbb{P}\left[y_{i}\in U_{f(i)} \ \text{for each} \ i\in [s]\right] \le \left(\frac{C'}{n}\right)^{s}.$$
\end{enumerate}
\end{lemma}


Note that the proof of Theorem~\ref{thm:bipartite perfect matching} from Lemma~\ref{lem:bipartite clustering} follows the same argument as the proof of Theorem~\ref{thm:hypergraph system spread} from Lemma~\ref{lem:clustering}, and is therefore omitted. 
Now it remains to prove Lemmas~\ref{lem:clustering} and~\ref{lem:bipartite clustering}. 
The proofs will be given in the next section.

\section{Proofs of random clustering lemmas}

\subsection{Concentration inequalities}

\begin{lemma}[Chernoff's inequality, \cite{Janson}]\label{lem: chernoff}
Let $X$ be a sum of independent Bernoulli random variables and $\lambda=\mathbb{E}(X)$.
Then for any $0<a<3/2$, we have
$$\mathbb{P}[X\ge (1+a)\lambda]\le e^{-a^2\lambda/3}~~~\text{and}~~~\mathbb{P}[X\le (1-a)\lambda]\le e^{-a^2\lambda/2}.$$
\end{lemma}


The following concentration inequality regarding random permutations is due to McDiarmid, which appears in~\cite[Chapter 16.3]{Michael-Bruce}. In this context, a \textit{choice} is defined to be either (a) the outcome of a trial
or (b) the position that a particular element gets mapped to in a permutation.

\begin{lemma}[McDiarmid's Inequality]\label{lem:McDiarmid's Inequality}
Let $X$ be a non-negative random variable, not identically 0, which is determined by $n$ independent trials $T_{1}, T_{2}, \dots , T_{n}$ and $m$ independent permutations $\Pi _{1}, \Pi _{2},\dots ,\Pi _{m}$ and satisfying the following for some $c, r>0$:
\begin{itemize}
    \item[1.] changing the outcome of any one trial can affect $X$ by at most $c$;
    \item[2.] interchanging two elements in any one permutation can affect $X$ by at most $c$;
    \item[3.] for any $s$, if $X\ge s$ then there is a set of at most $rs$ choices whose outcomes certify that $X\ge s$,
\end{itemize}
then for any $0\le t \le \mathbb{E}[X]$,
\begin{align*}
\mathbb{P}\left [|X-\mathbb{E}[X]|>t+60c\sqrt{r\mathbb{E}[X]} \right ]\le 4e^{-\frac{t^{2}}{8c^{2}r\mathbb{E}[X]}}.
\end{align*}
\end{lemma}

\begin{theorem}[Lemma 6.1 in \cite{Liebenau}]\label{lem:concentration tool}
Let $c>0$ and let $f$ be a function defined on the set of subsets of some set $U$ such that $|f(U_{1})-f(U_{2})|\le c$ whenever $|U_{1}|=|U_{2}|=m$ and $|U_{1}\cap U_{2}|=m-1$.
Let $A$ be a uniformly random $m$-subset of $U$.
Then for any $\alpha>0$, we have
\begin{align*}
\mathbb{P}\left[|f(A)-\mathbb{E}[f(A)]|\ge \alpha c \sqrt{m}\right]\le 2\exp(-2\alpha^{2}).
\end{align*}
\end{theorem}

\subsection{Random clustering in graph systems: Proof of Lemma~\ref{lem:clustering}}

In the proof of Lemma~\ref{lem:clustering}, we will use the following result in~\cite[Lemma 3.5]{Gupta} which is proved via a standard concentration argument.
\begin{lemma}[\cite{Gupta}]\label{lem:degree}
Let $k, \ell, d \in \mathbb{N}, 0 < \de'<\de < 1$ and $1/n, 1/\ell \ll 1/k, \delta - \delta'$. Let $H$ be an $n$-vertex $k$-graph with vertex set $V$ and suppose that $\deg(D, V) \geq \delta \binom{n-d}{k-d}$ for each $D \in \binom{V}{d}$. If $A \<  V$ is a vertex set of size $\ell$ chosen uniformly at random, then for every $D \in \binom{V}{d}$ we have
\[
\mathbb{P}\left[ \deg(D, A) < \delta' \binom{\ell-d}{k-d}\right] \leq 2 \exp\left(-\ell (\delta - \delta')^2 / 2\right).
\]
\end{lemma}

For the remainder of this section, choose $1/n\ll 1/C' \ll 1/C \ll \alpha,  1/k, 1/s, 1/t\ll 1$ as in Lemma~\ref{lem:clustering}. We have the following setup:
\begin{enumerate}
    \item  Let $\mathbf{H} =\{H_{1}, H_{2}, \dots ,H_{tn}\}$ be a $k$-graph system with the vertex set $V(\mathbf{H}):=\{v_{1}, v_{2}, \dots , v_{sn}\}$ and the color set $\mathcal{C}(\mathbf{H})=[tn]$ satisfying $\delta_{d}(\mathbf{H})\ge (\delta+\alpha)\binom{sn-d}{k-d}$ .
    \item  Recall that $r_{1}$ is $sC(C-1)$ plus the remainder when $sn$ is divided by $sC(C-1)$ and $r_{2}$ is $tC(C-1)$ plus the remainder when $tn$ is divided by $tC(C-1)$. 
    \item  Let $P_{1}:=[r_{1}]$, $Q_{1}=[r_{2}]$, and for $i\ge2$, let $P_{i}:=[r_{1}+s(C-1)(i-1)]\setminus  \bigcup_{j=1}^{i-1} P_{j}=[r_{1}+s(C-1)(i-2)+1, r_{1}+s(C-1)(i-1)]$ and $Q_{i}:=[r_{2}+t(C-1)(i-1)]\setminus  \bigcup_{j=1}^{i-1} Q_{j}=[r_{2}+t(C-1)(i-2)+1, r_{2}+t(C-1)(i-1)]$.
Note that $\{P_{1}, P_{2}, \dots , P_{m_{1}}\}$ is a partition of $[sn]$ with $m_{1}=(sn-r_{1})/s(C-1)+1$ and $\{Q_{1}, Q_{2}, \dots , Q_{m_{2}}\}$ is a partition of $[tn]$ with $m_{2}=(tn-r_{2})/t(C-1)+1$. As $tr_{1}=sr_{2}$, we may write $m':=m_{1}=m_{2}$.
    \item  Let $\pi _{1}$ be a uniformly random permutation of $[sn]$ and $\pi _{2}$ be a uniformly random permutation of $[tn]$.
Let $V_{i}:=\{v_{j}: \pi_{1}(j)\in P_{i}\}$ and $C_{i}:=\{j: \pi_{2}(j)\in Q_{i}\}$. Also let $\mathbf{H}_{i}:=\{H_{j}[V_{i}]\mid j\in C_{i}\}$ for every $i\in[m']$ and each $\mathbf{H}_{i}$ is called a \textit{sub-system} or a \textit{cluster}.
    \item For every $S\<  V(\mathbf{H})$ of size at most $s$ and every $T\<  \mathcal{C}(\mathbf{H})$ of size at most $t$, let $\mathbf{H}_{i}^{S\uplus T}=\{H_{j}[V_{i}\cup S]: j\in C_{i}\cup T\}$ for every $i\in [m']$.
\end{enumerate}

We shall first prove that for every $i\in [m']$, w.h.p.~the sub-system $\mathbf{H}_{i}$ inherits the minimum degree condition, that is, $\de_d(\mathbf{H}_{i})\ge (\de+\alpha/2)\binom{|V_i|-d}{k-d}$. Actually we can prove this in a stronger sense as follows.
\begin{lemma}\label{lem:main}
For every $S\<  V(\mathbf{H})$ of size at most $s$ and every $T\<  \mathcal{C}(\mathbf{H})$ of size at most $t$, and every $i\in [m']$, we have
\[
\mathbb{P}\left [\delta_{d}(\mathbf{H}_{i}^{S\uplus T})\ge (\delta+\alpha/2)\binom{|V_{i}\cup S|-d}{k-d}\right ]\ge 1- \exp(-sC\alpha^{2}/200).
\]
\end{lemma}

\begin{proof}
Let $S\<  V(\mathbf{H})$ of size at most $s$ and $T\<  \mathcal{C}(\mathbf{H})$ of size at most $t$.
Note that each $V_{i}$ is sampled from $V(\mathbf{H})$ uniformly at random for $i\in [m']$.
Fix $c\in C_{i}\cup T$, for every $D\<  V(\mathbf{H})$ with $|D|=d$, let $E_{D,c}$ be the event that $\deg_{H_{c}}(D,V_{i})\ge (\delta+2\alpha/3)\binom{|V_{i}|-d}{k-d}$.
Since $V_{i}$ is chosen uniformly at random, we have $\mathbb{P}[E_{D,c}]\ge 1-2\exp(-sC\alpha^{2}/100)$ by Lemma~\ref{lem:degree} applied with $\delta:=\delta+\alpha, \delta':=\delta+2\alpha/3$, $\ell=r_1$ when $i=1$ and $\ell=s(C-1)$ when $i\in [2,m']$. Since $1/C\ll \alpha, 1/s$, by the union bound, we have

\begin{align*}
&\mathbb{P}\left [\delta_{d}(\mathbf{H}_{i}^{S\uplus T})<(\delta+\alpha/2)\binom{|V_{i}\cup S|-d}{k-d}\right]\\
=& \mathbb{P}\left [\bigcup_{(D,c)}^{} \deg_{H_{c}}(D,V_{i}\cup S)< (\delta+\alpha/2)\binom{|V_{i}\cup S|-d}{k-d}\right]\\
\le&\sum_{(D,c)} \mathbb{P}\left [\deg_{H_{c}}(D,V_{i}\cup S)< (\delta+\alpha/2)\binom{|V_{i}\cup S|-d}{k-d}\right]\\
\le &\sum_{(D,c)} \mathbb{P}\left [\deg_{H_{c}}(D,V_{i})< (\delta+\alpha/2)\binom{|V_{i}\cup S|-d}{k-d}\right]\\
\le &\sum_{(D,c)} \mathbb{P}\left [\deg_{H_{c}}(D,V_{i})< (\delta+2\alpha/3)\binom{|V_{i}|-d}{k-d}\right]\\
\le &\binom{|V_{i}\cup S|}{d}\cdot |C_{i}\cup T|\cdot2\exp(-sC\alpha^{2}/100)\\
\le & \exp(-sC\alpha^{2}/200),
\end{align*}
as desired.
\end{proof}

The next lemma shows that w.h.p.~the sub-system $\mathbf{H}_{1}$ inherits the minimum degree condition, which can be immediately derived from Lemma~\ref{lem:main} applied with $S=T=\emptyset$ .
\begin{lemma}\label{lem:degree of V_{1}}
With probability at least $99/100$, $\delta_{d}(\mathbf{H}_{1})\ge (\delta+\alpha/2)\binom{|V_{1}|-d}{k-d}$.
\end{lemma}

In our random redistribution argument, we will assign every \textit{good} sub-system an $s$-set $S$ from the vertex set $V(\mathbf{H})$ and a $t$-set $T$ from the color set $\mathcal{C}(\mathbf{H})$ such that the minimum degree condition can be maintained (the definition of good and bad sub-system will be given later). The following lemma shows that w.h.p.~for most sub-systems $\mathbf{H}_{i}$, there are many choices of $(S,T)$ which are assignable to $\mathbf{H}_{i}$, that is, $\mathbf{H}_{i}^{S\uplus T}$ inherits the minimum degree condition.

\begin{lemma}\label{lem:vertex}
With probability at least $99/100$, for all but at most $\exp(-sC\alpha^{2}/1000)m'$ many $i\in [2,m']$, there are at least $(1-\exp(-sC\alpha^{2}/500))\binom{sn}{s}\binom{tn}{t}$ choices $(S,T)$ with $S\<  V(\mathbf{H})$ of size $s$ and $T\<  \mathcal{C}(\mathbf{H})$ of size $t$, such that
\begin{align*}
\delta_{d}(\mathbf{H}_{i}^{S\uplus T})\ge (\delta+\alpha/2)\binom{sC-d}{k-d}.
\end{align*}
\end{lemma}

\begin{proof}
By Lemma~\ref{lem:main}, for any fixed $i\in [2,m']$, $T\<  \mathcal{C}(\mathbf{H})$ of size $t$ and $S\<  V(\mathbf{H})$ of size $s$, we have that $\mathbb{P}\left [\delta_{d}(\mathbf{H}_{i}^{S\uplus T})< (\delta+\alpha/2)\binom{sC-d}{k-d}\right ]< \exp(-sC\alpha^{2}/200)$.
Call $S\cup T$ \textit{bad} for $\mathbf{H}_{i}$ if 
$\delta_{d}(\mathbf{H}_{i}^{S\uplus T})< (\delta+\alpha/2)\binom{sC-d}{k-d}$.
So for every $i\in [2,m']$, let $Y_i$ be the random variable counting the number of $(s+t)$-sets $S\cup T$ which are bad for $\mathbf{H}_{i}$. Then
\begin{align*}
\mathbb{E}[Y_i]\le \exp (-sC\alpha^{2}/200)\binom{sn}{s}\binom{tn}{t}.
\end{align*}
By Markov's inequality, for any $i\in [2,m']$, we have
\begin{align*}
\mathbb{P}\left[Y_i\ge  \exp (-sC\alpha^{2}/500)\binom{sn}{s}\binom{tn}{t} \right ]\le \exp (-sC\alpha^{2}/500).
\end{align*}
Let $X$ denote the number of $\mathbf{H}_{i}$ such that $Y_i\ge  \exp (-sC\alpha^{2}/500)\binom{sn}{s}\binom{tn}{t}$.
Then by linearity of expectation,
\begin{align*}
\mathbb{E}[X]\le \exp (-sC\alpha^{2}/500)m'.
\end{align*}
By Markov's inequality, we have
$
\mathbb{P}[X\ge \exp (-sC\alpha^{2}/1000)m']\le \exp (-sC\alpha^{2}/1000),
$
implying the desired statement.

\end{proof}



The following lemma tells that w.h.p.~for every choice of $(S,T)$ with $S\<  V(\mathbf{H})$ of size  $s$ and 
$T\<  \mathcal{C}(\mathbf{H})$ of size $t$, there are many sub-systems $\mathbf{H}_{i}$ such that $\mathbf{H}_{i}^{S\uplus T}$ inherits the minimum degree condition.
\begin{lemma}\label{lem:subhypergraph system}
With probability at least $99/100$, for every $S\< V(\mathbf{H})$ of size at most $s$ and $T\< \mathcal{C}(\mathbf{H})$ of size at most $t$, there are at least $(1-1/C^{2})m'$ indices $i\in [2,m']$ for which 
\begin{align*}
\delta_{d}(\mathbf{H}_{i}^{S\uplus T})\ge (\delta+\alpha/2)\binom{|V_i\cup S|-d}{k-d}.
\end{align*}

\end{lemma}

\begin{proof}
For every $S\< V(\mathbf{H})$ of size at most $s$ and every $T\< \mathcal{C}(\mathbf{H})$ of size at most $t$, let $Y_{S,T}$ denote the number of sub-systems $\mathbf{H}_{i}$ such that $\delta_{d}(\mathbf{H}_{i}^{S\uplus T})\ge (\delta+\alpha/2)\binom{|V_i\cup S|-d}{k-d}$.
By Lemma~\ref{lem:main}, we have
\begin{align*}
\mathbb{E}\left [Y_{S,T} \right ]\ge (1-\exp(-sC\alpha^{2}/200))m'.
\end{align*}
Note that interchanging two elements of the random permutation $\pi_{1}$ (or $\pi_{2}$) either has no effect on $V_{1}, V_{2}, \dots , V_{m'}$ or changes two random sets $V_{i}$ and $V_{j}$ and therefore affects the value of $Y_{S,T}$ by at most 2.
Also, if $Y_{S,T}\ge s$, 
then there is a set of at most $s^2C^{2}$ choices whose outcomes certify that $Y_{S,T}\ge s$.
Therefore, by Lemma~\ref{lem:McDiarmid's Inequality} applied with $c=2, r=sC^{2}$, and $t=m'/C^{3}$, and the fact $m'>n/C^{2}$, we have
\begin{align*}
\mathbb{P}\left [Y_{S,T}<(1-1/C^{2})m' \right ]\le 4\exp\left (-\frac{(m'/C^{3})^{2}}{32sC^{2}m'}\right)\le \exp\left(-n/C^{12} \right).
\end{align*}
By the union bound, we have that with probability at least $99/100$, $Y_{S,T}\ge (1-1/C^{2})m'$ for every $S\< V(\mathbf{H})$ of size at most $s$ and $T\< \mathcal{C}(\mathbf{H})$ of size at most $t$, as desired.
\end{proof}

Now we are ready to decode the proof of Lemma~\ref{lem:clustering}.

\begin{proof}[Proof of Lemma~\ref{lem:clustering}]
Let $E_{1}, E_2, E_3$ be the events for the respective properties in Lemmas~\ref{lem:degree of V_{1}},~\ref{lem:vertex} and~\ref{lem:subhypergraph system}. 
We will condition on $E_{1}\cap E_{2}\cap E_{3}$, which holds with probability at least 97/100.
Recall that $\pi_1$ is a uniformly random permutation of $[sn]$ and $\pi_2$ is a uniformly random permutation of $[tn]$. Now we only consider $\pi_1,\pi_2$ such that $(\pi_1,\pi_2)\in E_{1}\cap E_{2}\cap E_{3}.$
Recall again that $\{V_1,\ldots,V_{m'}\}$ is the partition of vertex set $V(\mathbf{H})$, $\{C_1,\ldots,C_{m'}\}$ is the partition of color set $\mathcal{C}(\mathbf{H})$,
and $\mathbf{H}_{i}=\{H_{j}[V_{i}] \mid j\in C_{i}\}$ (or $\mathbf{H}_{i}=(V_i,C_i)$ for brevity).

Define a family $\mathcal{F}_{\pi_{1},\pi_{2}}$ of \textit{bad} sub-systems $\mathbf{H}_{i}$, where $i\in [2,m']$, for which at least one of the following properties holds:
\begin{itemize}
\item[(i)] $\delta_{d}(\mathbf{H}_{i})<(\delta+\alpha/2)\binom{|V_{i}|-d}{k-d}$;
\item[(ii)] $\delta_{d}(\mathbf{H}_{i}^{S\uplus T})<(\delta+\alpha/2)\binom{sC-d}{k-d}$ for at least $\exp(-sC\alpha^{2}/500)\binom{sn}{s}\binom{tn}{t}$ choices $(S,T)$ with $S\< V(\mathbf{H})$ of size $s$ and $T\< \mathcal{C}(\mathbf{H})$ of size $t$.
\end{itemize}
Let $m'_{1}=|V_{2}\cup\dots \cup V_{m'}|/sC$ and $m'_{2}=|C_{2}\cup\dots \cup C_{m'}|/tC$. 
Let $m-1:=m'_{1}=m'_{2}=(m'-1)(1-\tfrac{1}{C})$ as $m'_{1}=m'_{2}$ is a positive integer by the choice of $r_{1}$ and $r_{2}$.
We claim that if $(\pi_{1}, \pi_{2})\in E_{1}\cap E_{2}\cap E_{3}$, then
we can add extra elements among $\mathbf{H}_{2}, \dots , \mathbf{H}_{m'}$ arbitrarily to $\mathcal{F}_{\pi_{1},\pi_{2}}$ to make sure $|\mathcal{F}_{\pi_{1},\pi_{2}}|= m'-1-(m-1)=m'-m$.
Indeed, there are at most $m'/C^{2}$ bad sub-systems of type (i) by Lemma~\ref{lem:subhypergraph system} applied with $S=T=\emptyset $, and at most $\exp(-sC\alpha^{2}/1000)m'$ bad sub-systems of type (ii) by Lemma~\ref{lem:vertex}.

We define a balanced $(s+t+1)$-partite $(s+t+1)$-graph $H_{\pi_{1},\pi_{2}}$ with parts $\{\hat{\mathcal{F}}, V_1',\cdots, V_s', C_1',\cdots, C_t'\}$ satisfying the following conditions:
\begin{enumerate}
    \item  $\hat{\mathcal{F}}=\{\mathbf{H}_{2},\dots , \mathbf{H}_{m'}\} \setminus  \mathcal{F}_{\pi_{1},\pi_{2}}$;
    \item  $\{V_{1}',\dots , V_{s}'\}$ and $\{C_{1}',\dots , C_{t}'\}$ are arbitrary (balanced) partitions of the vertex set $\bigcup_{\mathbf{H}_i\in \mathcal{F}_{\pi_{1},\pi_{2}}}V_i$ and  color set $\bigcup_{\mathbf{H}_i\in \mathcal{F}_{\pi_{1},\pi_{2}}}C_i$, respectively;
    \item each part has exactly $m-1$ elements;
    \item for all $S=\{v_{i}\in V_{i}': i\in [s]\}$, $T=\{c_{j}\in C_{j}' : j\in [t]\}$ and $\mathbf{H}_{i}\in \hat{\mathcal{F}}$,  we put $S\cup T \cup \{\mathbf{H}_{i}\}$ into $E(H_{\pi_{1},\pi_{2}})$ whenever $\delta_{d}(\mathbf{H}_{i}^{S\uplus T})\ge (\delta+\alpha/2)\binom{sC-d}{k-d}$.
\end{enumerate}
The first three items are indeed valid as $|\hat{\mathcal{F}}|=m'-1-|\mathcal{F}_{\pi_{1},\pi_{2}}|=m-1$. 
Without loss of generality assume that $\hat{\mathcal{F}}=\{\mathbf{H}_{2},\dots , \mathbf{H}_{m}\}$.
Also, 
as $|\bigcup_{H_i\in \mathcal{F}_{\pi_{1},\pi_{2}}}V_i|=(m'-1-(m-1))\cdot s(C-1)=s(m-1)$, we can divide it arbitrarily into $s$ parts $V_{1}',\dots , V_{s}'$, each of size $m-1$, and similarly for $C_{1}',\dots , C_{t}'$.

Conditioning on $E_{1}\cap E_{2}\cap E_{3}$, we claim that $\delta_{\{1\}}(H_{\pi_{1},\pi_{2}})\cdot (m-1)+\de_{\{2,\cdots, s+t+1\}}(H_{\pi_{1},\pi_{2}})\cdot (m-1)^{s+t}\ge (1+\alpha)(m-1)^{s+t+1}$, which leads to the degree condition of Theorem~\ref{thm:multipartite spread}. Indeed, for every $\mathbf{H}_{i}\in \hat{\mathcal{F}}$, by Lemma~\ref{lem:vertex}, there are at most $\exp(-sC\alpha^{2}/500)\binom{sn}{s}\binom{tn}{t}$ choices $(S,T)$ with $S\< V(\mathbf{H})$ of size $s$ and $T\< \mathcal{C}(\mathbf{H})$ of size $t$ such that $\delta_{d}(\mathbf{H}_{i}^{S\uplus T})<(\delta+\alpha/2)\binom{sC-d}{k-d}$. 
So by the choice $1/C \ll \alpha, 1/s,1/t$, we obtain 
\begin{align*}
\de_{\{1\}}(H_{\pi_{1},\pi_{2}})\ge (m-1)^{s+t}-\exp(-sC\alpha^{2}/500)\binom{sn}{s}\binom{tn}{t}\ge (1/2+\alpha/2)(m-1)^{s+t}.
\end{align*}
For every vertex subset $S$ of size $s$ from different parts of $V_1',\cdots, V_s'$ and color subset $T$ of size $t$ from different parts of $C_1',\cdots, C_t'$, Lemma~\ref{lem:subhypergraph system} implies there are at least $(1-1/C^2)m'$ indices $i\in [2,m']$ such that $\delta_{d}(\mathbf{H}_{i}^{S\uplus T})\ge (\delta+\alpha/2)\binom{sC-d}{k-d}$. As $1/C \ll \alpha$, we have
\begin{align*}
\de_{\{2,\cdots, s+t+1\}}(H_{\pi_{1},\pi_{2}}) \ge (1-1/C^{2})m'-|\mathcal{F}_{\pi_{1},\pi_{2}}|\ge (1/2+\alpha/2)(m-1).
\end{align*}
Therefore, by Theorem~\ref{thm:multipartite spread}, there is a $(C_{\ref{thm:multipartite spread}}/(m-1))$-vertex-spread distribution on embeddings of perfect matchings in $H_{\pi_{1},\pi_{2}}$.

We are ready to define the  vertex sets $U_{i}$ and color sets $W_{i}$ as desired in Lemma~\ref{lem:clustering}.
First, randomly choose $(\pi_1,\pi_2)\in E_1\cap E_2 \cap E_3$.
Then, sample $M_{\pi_{1},\pi_{2}}$ from the $(C_{\ref{thm:multipartite spread}}/(m-1))$-vertex-spread distribution on embeddings of perfect matchings of $H_{\pi_{1},\pi_{2}}$.
For every $\mathbf{H}_{i}=(V_{i}, C_{i})\in \hat{\mathcal{F}}$, let $S_i=\{v_{i,1}, \dots , v_{i,s}\}$ and $T_i=\{c_{i,1}, \dots , c_{i,t}\}$ be the vertex and color sets that together with $\mathbf{H}_{i}$ form an edge in $M_{\pi_{1},\pi_{2}}$ where $v_{i,j}\in V_j'$ for $j\in [s]$ and $c_{i,j}\in C_j'$ for 
 $j\in [t]$.
Let $U_{i}=V_{i}\cup S_i$ , $W_{i}=C_{i}\cup T_i$ for $i\in \{2,\dots,m\}$ and let $U_{1}=V_{1}$, $W_{1}=C_{1}$.
In this case, it remains to verify that the resulting partitions $\mathcal{U}:=\{U_{1}, U_{2}, \dots, U_{m}\}$ and $\mathcal{W}:=\{W_{1}, W_{2}, \dots, W_{m}\}$ for the vertex set $V(\mathbf{H})$ and color set $\mathcal{C}(\mathbf{H})$, respectively, satisfy the properties~\ref{size}\ref{degree}\ref{location} in Lemma~\ref{lem:clustering}.
As~\ref{size} and~\ref{degree} easily follow from our construction, we only verify~\ref{location}.
Let $y_{1}, y_{2},\dots , y_{r}\in V(\mathbf{H})$, $c_{1}, c_{2},\dots , c_{\ell }\in \mathcal{C}(\mathbf{H})$ and let $f: [r]\to [m]$, $g: [\ell]\to [m]$ be arbitrarily given as in~\ref{location}.
For any $i\in [r], j\in [\ell]$, we denote by $D_{i,j}$ the event that $y_{i}\in U_{f(i)}=V_{f(i)}\cup S_{f(i)}$ and $c_{j }\in W_{g(j)}=C_{g(j)}\cup T_{g(j)}$, $D_{i,j}^{1}$ the event that $y_{i}\in V_{f(i)}$ and $c_{j }\in C_{g(j)}$, and $D_{i,j}^{2}$ the event that $y_{i}\in S_{f(i)}$ and $c_{j }\in T_{g(j)}$. 
Note that the event $D_{i,j}^{1}$ is determined only by the permutations $\pi_{1}$ and $\pi_{2}$, and the event $D_{i,j}^{2}$ happens only when  there are $(V_{f(i)},C_{f(i)}), (V_{g(j)},C_{g(j)})\in \hat{\mathcal{F}}$ such that $\{y_{i},(V_{f(i)},C_{f(i)})\}$ and $\{c_{j},(V_{g(j)},C_{g(j)})\}$ are in edges of $M_{\pi_{1},\pi_{2}}$.
Thus, for every $S_{1}\< [r]$ and $S_{2}\< [\ell]$, we have
\begin{align*}
\mathbb{P}\left[\bigcap_{\substack{i \in S_1, j \in S_2}} D_{i,j}^1\right] 
&\le \frac{(sC-s)^{|S_{1}|} (tC-t)^{|S_{2}|} (sn-|S_{1}|)! (tn-|S_{2}|)!}{(sn)!\cdot (tn)!} \cdot \frac{1}{\mathbb{P}[E_{1}\cap E_{2}\cap E_{3}]}\le \frac{100}{97}\left(\frac{C'}{n}\right)^{|S_{1}|+|S_{2}|}.
\end{align*}
and for every $(\pi_{1}, \pi_{2})\in E_{1}\cap E_{2}\cap E_{3}$,
\begin{align*}
\mathbb{P}\left[\bigcap_{\substack{i \in [r] \setminus S_1 \\ j \in [\ell] \setminus S_2}} D_{i,j}^{2} \bigg | (\pi_{1}, \pi_{2}) \right ] 
&\le \left(\frac{C_{\ref{thm:multipartite spread}}}{m-1}\right)^{r-|S_{1}|}\cdot \left(\frac{C_{\ref{thm:multipartite spread}}}{m-1}\right)^{\ell-|S_{2}|}
\le \left(\frac{C'}{n}\right)^{r -|S_{1}|+\ell-|S_{2}|}.
\end{align*}
Therefore, for every $S_{1}\< [r]$ and $S_{2}\< [\ell]$, we have
\begin{align*}
\mathbb{P}\left[\bigcap_{\substack{i \in  S_1 \\ j \in  S_2}} D_{i,j}^{1} \cap  \bigcap_{\substack{i \in [r] \setminus S_1 \\ j \in [\ell] \setminus S_2}} D_{i,j}^{2} \right ] 
&=  \mathbb{P}\left[\bigcap_{\substack{i \in  S_1 \\ j \in  S_2}} D_{i,j}^{1}  \right ] \cdot \mathbb{P}\left[\bigcap_{\substack{i \in [r] \setminus S_1 \\ j \in [\ell] \setminus S_2}} D_{i,j}^{2} \bigg | \bigcap_{\substack{i \in  S_1 \\ j \in  S_2}} D_{i,j}^{1} \right ] \\ 
&\le \frac{100}{97}\left(\frac{C'}{n}\right)^{|S_{1}|+|S_{2}|} \cdot \left(\frac{C'}{n}\right)^{r -|S_{1}|+\ell-|S_{2}|}\\
&\le \left(\frac{C'}{2n}\right)^{r+\ell }.
\end{align*}
Since $\mathbb{P}\big{[}y_{i}\in U_{f(i)}~\text{for all}~i\in [r] ~\text{and}~ c_{j}\in W_{g(j)}~\text{for all}~j\in [\ell]\big{]}=\mathbb{P}\big{[}\bigcap_{\substack{i \in  [r], j \in  [\ell]}} D_{i,j} {]}$ and 
\begin{align*}
\bigcap_{\substack{i \in  [r], j \in  [\ell]}} D_{i,j}
=\bigcup_{\substack{ S_{1}\< [r]\\ S_{2}\< [\ell ]} }\left(\bigcap_{\substack{i \in  S_1 \\ j \in  S_2}} D_{i,j}^{1} \cap  \bigcap_{\substack{i \in [r] \setminus S_1 \\ j \in [\ell] \setminus S_2}}  D_{i,j}^{2}\right),
\end{align*}
the result follows by the union bound over the $2^{r+\ell }$ choices of $S_{1}\< [r]$ and $S_{2}\< [\ell]$, that is, 
\begin{align*}
\mathbb{P}\left[\bigcap_{\substack{i \in  [r], j \in  [\ell]}} D_{i,j} \right] \
\le\sum_{\substack{S_{1}\< [r]\\S_{2}\< [\ell]}} \mathbb{P}\left[\bigcap_{\substack{i \in  S_1 \\ j \in  S_2}} D_{i,j}^{1} \cap  \bigcap_{\substack{i\in [r]\setminus S_{1}\\j\in [\ell ]\setminus S_{2}}} D_{i,j}^{2}\right]
\le  2^{r+\ell }\cdot \left(\frac{C'}{2n}\right)^{r+\ell }=\left(\frac{C'}{n}\right)^{r+\ell}.
\end{align*}

\end{proof}


\subsection{Random clustering in balanced bipartite graphs: Proof of Lemma~\ref{lem:bipartite clustering}}

Fix $1/n \ll 1/C' \ll 1/C \ll \varepsilon $.
Let $r$ be $(C-1)C$ plus the remainder when $n$ is divided by $(C-1)C$.
Let $W_{1}^{1}=W_{1}^{2}:=[r]$ and 
$W_{i}^{p}:=[r+(C-1)(i-1)]\setminus  {\textstyle \bigcup_{j=1}^{i-1}} W_{j}^{p}$ for $i\ge 2$ and $p=1,2$.
Note that $\{W_{1}^{p}, \dots , W_{m'}^{p}\}$ is a partition of $[n]$ for $m'=(n-r)/(C-1)+1$.
Let $G$ be a balanced bipartite graph with vertex set $V(G)=V_1\cup V_2=\{v_{1},\dots, v_{n}\}\cup \{u_{1},\dots, u_{n}\}$ satisfying $\delta_{\{1\}}(G)+\delta_{\{2\}}(G)\ge (1+\varepsilon )n$.
Let $\pi_{1}$, $\pi_{2}$ be permutations of $[n]$ chosen uniformly and independently, and write $D_{i}:=\{v_{j}: \pi_1(j)\in W_{i}^{1}\}\cup \{u_{j}: \pi_2(j)\in W_{i}^{2}\}$ for every $i\in[m']$. Note that $|D_1|=2r$ and $|D_i|=2(C-1)$ for each $2\le i\le m'$.\medskip

We first prove that every cluster $D_i$ inherits the degree condition w.h.p.~in a strong sense as follows.
\begin{lemma}\label{lem:degree condition bipartite}
For every $T\< V(G)$ with $|T\cap V_{j}|\le 1$ for each $j\in [2]$, and every $i\in [m']$, 
\begin{align*}
\mathbb{P}\left[\delta_{\{1\}}(G[D_{i}\cup T])+\delta_{\{2\}}(G[D_{i}\cup T])< (1+\varepsilon/2)\frac{|D_i|}{2}\right]< \exp(-C\varepsilon^{2}/8).
\end{align*}
\end{lemma}

To prove Lemma~\ref{lem:degree condition bipartite}, 
we need the following result obtained by a standard concentration inequality.

\begin{lemma}\label{lem: local degree condition}
Let $\ell \in \mathbb{N}$, $0<\varepsilon <1$ and $1/n, 1/\ell \ll \eps$.  
Let $G=(V_{1}, V_{2})$ be a bipartite graph with $|V_{1}|=|V_{2}|=n$ and suppose that $\delta_{\{1\}}(G)+\delta_{\{2\}}(G)\ge (1+\varepsilon )n$.
If $A\< V(G)$ is a vertex subset with $|A\cap V_{i}|=\ell $ for each $i\in [2]$ chosen uniformly at random, then for every $v\in V_{1}$ and $u\in V_{2}$ we have
\begin{align*}
\mathbb{P}\left[\deg(v, A)+ \deg(u, A)<(1 +\tfrac{\varepsilon}{2})\ell \right]\le 2\exp(-\ell \varepsilon^{2}/4).
\end{align*}
\end{lemma}


\begin{proof}[Proof of Lemma~\ref{lem: local degree condition}]
Given $\ell\in \mathbb{N}$.
Let $\mathcal{A}$ be the family of all $2\ell$-subsets $U\<V(G)$ with $|U\cap V_i|=\ell$ for $i\in[2]$.
For every $v\in V_{1}$ and $u\in V_{2}$, let $f: \mathcal{A} \rightarrow \mathbb{R}$ be defined by $f(U)=\deg(v,U)+\deg(u,U)$ for each $U\in \mathcal{A}$.
Observe that $|f(U)-f(U')|\le 1$ for any $U$, $U'\in \mathcal{A}$ with $|U\cap U'|=2\ell -1$.
Note that for any vertex $v\in V_{1}$ (similarly for $u\in V_{2}$), the probability that $v\in A$ is at least $\frac{\binom{n-1}{\ell-1}}{\binom{n}{\ell}}=\frac{\ell}{n}$.
Thus by linearity of expectation every $v\in V_{1}$ has
\begin{align*}
\mathbb{E}\left[\deg (v,A)\right]=\sum_{u\in N(v)}\mathbb{P}[u\in A]\ge\deg(v,V(G))\cdot \tfrac{\ell}{n}\ge  \delta_{\{1\}}(G)\cdot \tfrac{\ell}{n} .
\end{align*}
Similarly, for each $u\in V_{2}$, 
\begin{align*}
\mathbb{E}\left[\deg (u,A)\right]=\sum_{v\in N(u)}\mathbb{P}[v\in A]\ge\deg(u,V(G))\cdot \tfrac{\ell}{n}\ge \delta_{\{2\}}(G)\cdot \tfrac{\ell}{n} .
\end{align*}
Therefore,  
$\mathbb{E}[f(A)]
=\mathbb{E}\left[\deg (v,A)\right]+\mathbb{E}\left[\deg (u,A)\right]
\ge \left(\delta_{\{1\}}(G) +\delta_{\{2\}}(G)\right)\cdot \tfrac{\ell}{n} 
\ge (1+\varepsilon )\ell$. By applying Theorem~\ref{lem:concentration tool} with $c=1$, $m=2\ell$, $\alpha=\frac{\varepsilon}{2} \sqrt{\frac{\ell}{2}}$, we obtain that
\begin{align*}
\mathbb{P}\left[f(A)<(1 +\tfrac{\varepsilon }{2})\ell \right]
&=\mathbb{P}[f(A)<(1+\varepsilon)\ell-\alpha c\sqrt{m}]\\
&\le \mathbb{P}[f(A)\le\mathbb{E}[f(A)]-\alpha c\sqrt{m}]\\
&\le \mathbb{P}\left[|f(A)-\mathbb{E}[f(A)]|\ge \alpha c \sqrt{m}\right]\\
&\le 2\exp(-\ell \varepsilon^{2}/4),
\end{align*}
as desired.
\end{proof}

\begin{proof}[Proof of Lemma~\ref{lem:degree condition bipartite}]
For every $v\in V_{1}\cap(D_{i}\cup T)$ and $u\in V_{2}\cap(D_{i}\cup T)$, we have
\begin{align*}
\mathbb{P}\left[\deg(v, D_{i}\cup T)+\deg(u, D_{i}\cup T)<(1+\tfrac{\varepsilon }{2})\tfrac{|D_i|}{2}\right]\le \mathbb{P}\left[\deg(v, D_{i})+\deg(u, D_{i})<(1+\tfrac{\varepsilon}{2})\tfrac{|D_i|}{2}\right].
\end{align*}
Then, by the union bound and Lemma~\ref{lem: local degree condition}, we have
\begin{align*}
&\mathbb{P}\left[\delta_{\{1\}}(G[D_{i}\cup T])+\delta_{\{2\}}(G[D_{i}\cup T])< (1+\tfrac{\varepsilon}{2} )\tfrac{|D_i|}{2}\right]\\
\le& |D_i\cup T|^2\cdot \mathbb{P}\left[\deg(v, D_{i}\cup T)+\deg(u, D_{i}\cup T)<(1+\tfrac{\varepsilon }{2})\tfrac{|D_i|}{2}\right]\\
\le& |D_i\cup T|^2\cdot 2\exp(-C\varepsilon^{2}/4)\\
\le& \exp(-C\varepsilon^{2}/8),
\end{align*}
as desired.
\end{proof}

The next lemma ensures that w.h.p.~$D_1$ inherits the minimum degree condition, which can be deduced from Lemma~\ref{lem:degree condition bipartite} applied with $T=\emptyset $.
\begin{lemma}\label{lem:degree condition bipartite1}
With probability at least $99/100$, 
$\delta_{\{1\}}(G[D_{1}])+\delta_{\{2\}}(G[D_{1}])\ge (1+\tfrac{\varepsilon}{2})\tfrac{|D_1|}{2}$.
\end{lemma}

In our random redistribution argument, we will assign every ``good" random cluster $D_i$ two vertices, one from $V_1$ and one from $V_2$ such that the minimum degree condition can be retained. The following lemma shows that w.h.p.~for every two vertices from distinct parts, there are many choices of $D_i$ which are assignable.

\begin{lemma}\label{lem:auxiliary bipartite degree 1}
With probability at least $99/100,$ for every $T\< V(G)$ with $|T\cap V_{j}|\le 1$ for each $j\in [2]$, there are at least $(1-1/C)m'$ clusters $D_{i}$, where $i\in [2,m']$, for which 
\begin{align*}
\delta_{\{1\}}(G[D_{i}\cup T])+\delta_{\{2\}}(G[D_{i}\cup T])\ge (1+\tfrac{\varepsilon}{2} )C.
\end{align*}
\end{lemma}
\begin{proof}
For every $T\< V(G)$ with $|T\cap V_{j}|\le 1$ for each $j\in [2]$, let $X_{T}$ denote the number of random clusters $D_{i}$ such that the degree condition above holds.
By the linearity of expectation and Lemma~\ref{lem:degree condition bipartite}, we have
\begin{align*}
\mathbb{E}[X_{T}]\ge (1-\exp(-C\varepsilon^{2}/8))m'.
\end{align*}
Note that interchanging two elements of the random permutation $\pi_{1}$ (or $\pi_{2}$) either has no effect on $D_{2},\dots ,D_{m'}$ or changes two random clusters $D_{i}$ and $D_{j}$, therefore it can effect the value of $X_{T}$ by at most 2.
Also, if $X_{T}\ge s$, 
then there is a set of at most $2Cs$ choices whose outcomes certify that $X_{T}\ge s$.
Therefore, by Lemma~\ref{lem:McDiarmid's Inequality} applied with $c=2, r=2C$, and $t=m'/C^{3}$, and the fact $m'>n/C^2$, we have
\begin{align*}
\mathbb{P}\left [X_{T}<(1-1/C)m' \right ]
&\le \mathbb{P}[X_{T}\le \mathbb{E}[X_{T}]-t-60c\sqrt{r\mathbb{E}[X_{T}]}]\\
&\le 4\exp\left (-\frac{(m'/C^{3})^{2}}{64Cm'}\right)\\
&\le \exp\left(-n/C^{12} \right).
\end{align*}
By the union bound, with the probability at least 99/100 we have that $X_{T}\ge (1-1/C)m'$ for  every $T\< V(G)$ with $|T\cap V_{j}|\le 1$ for each $j\in [2]$.
\end{proof}

The following lemma tells that w.h.p.~for most clusters $D_i$, there are many choices of vertices which are assignable for $D_i$.

\begin{lemma}\label{lem:auxiliary bipartite degree 2}
With probability at least $99/100$, for all but at most $\exp(-C\varepsilon^{2}/32)m'$ many $i\in [2,m']$, there are at least $(1-2\exp(-C\varepsilon^{2}/16))n^{2}$ subsets $T\< V(G)$ with $|T\cap V_{j}|\le 1$ for each $j\in [2]$ such that
\begin{align*}
\delta_{\{1\}}(G[D_{i}\cup T])+\delta_{\{2\}}(G[D_{i}\cup T])\ge (1+\tfrac{\varepsilon}{2} )C.
\end{align*}
\end{lemma}

\begin{proof}
Let $\mathcal{T}=\{T\subseteq V(G): |T\cap V_{j}|\le 1~\text{for each}~j\in [2]\}$.
For any $T\in \mathcal{T}$ and $i\in [2,m']$, 
by Lemma~\ref{lem:degree condition bipartite} we have that 
\begin{align*}
\mathbb{P}\left[\delta_{\{1\}}(G[D_{i}\cup T])+\delta_{\{2\}}(G[D_{i}\cup T])< (1+\tfrac{\varepsilon}{2} )C\right]< \exp(-C\varepsilon^{2}/8).
\end{align*}
We call $T$ \textit{bad} for $D_{i}$ if $\delta_{\{1\}}(G[D_{i}\cup T])+\delta_{\{2\}}(G[D_{i}\cup T])< (1+\tfrac{\varepsilon}{2} )C$.
By the linearity of expectation, we have
\begin{align*}
\mathbb{E}\left[|\{T\in \mathcal{T}: T \ \text{is bad for}\ D_{i}\}|\right]\le &(n^2+2n)\cdot \exp(-C\varepsilon^{2}/8)\\
\le&2n^2\cdot \exp(-C\varepsilon^{2}/8)
\end{align*}
for any $i\in[2,m']$. By Markov's inequality, we have
\begin{align*}
\mathbb{P}\left[|\{T\in \mathcal{T}: T \ \text{is bad for}\ D_{i}\}|\ge 2n^{2}\cdot \exp(-C\varepsilon^{2}/16)\right]\le \exp(-C\varepsilon^{2}/16).
\end{align*}
Let $Y$ denote the number of random clusters $D_{i}$ such that $|\{T\in \mathcal{T}: T \ \text{is bad for}\ D_{i}\}|\ge 2n^{2}\cdot \exp(-C\varepsilon^{2}/16)$.
Then $\mathbb{E}[Y]\le \exp(-C\varepsilon^{2}/16)m'$.
Again by Markov's inequality, 
we obtain \[\mathbb{P}[Y\ge \exp(-C\varepsilon^{2}/32)m']\le \exp (-C\varepsilon^{2}/32),\]
yielding the desired statement.
\end{proof}

Now we are ready to prove Lemma~\ref{lem:bipartite clustering}.

\begin{proof}[Proof of Lemma~\ref{lem:bipartite clustering}]
Let $E_{1}, E_{2}, E_{3}$ be the events that the properties in Lemmas~\ref{lem:degree condition bipartite1},~\ref{lem:auxiliary bipartite degree 1} and~\ref{lem:auxiliary bipartite degree 2} hold respectively.
We will condition on $E_{1}\cap E_{2} \cap E_{3}$, which holds with the probability at least $97/100$.

Sample two permutations $\pi_{1}$, $\pi_{2}$ independently from the uniform distribution on permutation of $[n]$ conditional on $E_{1}\cap E_{2} \cap E_{3}$. Write $\pi:=(\pi_{1},\pi_{2})$
and define a family of \textit{bad} random clusters $\mathcal{D}_{\pi}\< \{D_{2}, \dots , D_{m'}\}$ by including each $D_{i}$ for which we have any of the following properties:
\begin{itemize}
\item [(1)] $\delta_{\{1\}}(G[D_{i}])+\delta_{\{2\}}(G[D_{i}])<(1+\frac{\varepsilon}{2})(C-1)$,
\item [(2)] $\delta_{\{1\}}(G[D_{i}\cup T])+\delta_{\{2\}}(G[D_{i}\cup T])<(1+\frac{\varepsilon}{2})C$ for at least $2\exp(-C\varepsilon^{2}/16)n^{2}$ subsets $T\< V(G)$ with $|T\cap V_{i}|\le 1$ for each $i\in [2]$.
\end{itemize}
Let $m-1=|D_{2}\cup\dots\cup D_{m'}|/(2C)$ and note that $m$ is a positive integer by the choice of $r$.
If necessary, we can add extra good clusters of $\{D_{2}, \dots , D_{m'}\}$ arbitrarily to $\mathcal{D}_{\pi}$ to make sure $\mathcal{D}_{\pi}$ has size exactly $m'-1-(m-1)=m'-m$.
Since we condition on $E_{1}\cap E_{2} \cap E_{3}$, there are very few bad clusters.
Indeed, by Lemma~\ref{lem:auxiliary bipartite degree 1} applied with $T=\emptyset$, there are at most $m'/C$ clusters satisfying  condition (1) 
and by Lemma~\ref{lem:auxiliary bipartite degree 2}, there are at most $\exp(-C\varepsilon^{2}/32)m'$ clusters satisfying condition (2). 
Therefore the family $\mathcal{D}_{\pi}$ is of size exactly $m'-m$.\medskip

Next we iteratively redistribute the vertices of bad clusters into good ones in two rounds. We first define a bipartite graph $H^1_{\pi}$ between $A:={\textstyle \bigcup_{D_{i}\in \mathcal{D}_{\pi}}^{}}D_{i}\cap V_{1} $ and $B:=\{D_{2}, \dots , D_{m'}\}\setminus \mathcal{D}_{\pi}$, such that there is an edge between $v^{1}\in A$ and $D_{i}\in B$ if $\delta_{\{1\}}(G[D_{i}\cup \{v^{1}\}])+\delta_{\{2\}}(G[D_{i}\cup \{v^{1}\}])\ge(1+\tfrac{\varepsilon}{3})C$ and there are at least $(1-2\exp(-C\varepsilon^{2}/16))n$ vertices $v^{2}\in V_{2}$ such that $\delta_{\{1\}}(G[D_{i}\cup T])+\delta_{\{2\}}(G[D_{i}\cup T])\ge(1+\tfrac{\varepsilon}{3})C$ where $T=\{v^{1},v^{2}\}$. Note that $A$ is a set of vertices from bad clusters with $|A|=(m'-m)(C-1)=(m'-1)(C-1)-(m-1)C+(m-1)=m-1$ and $B$ consists of $m-1$ good clusters. 
By applying Lemma~\ref{lem:auxiliary bipartite degree 1} with $T=\{v^1\}$ and condition (2),
we have 
\begin{align*}
d_{H^1_{\pi}}(v^{1})\ge (1-1/C)m'-|\mathcal{D}_{\pi}|=(1-1/C)m'-(m'-m)\ge \tfrac{3}{4}(m-1)
\end{align*}
for each $v^1\in A$. 
By condition (2), each $D_i\in B$ has 
\begin{align*}
d_{H^1_{\pi}}(D_{i})\ge (m-1)-2\exp(-C\varepsilon^{2}/16)n\ge \tfrac{3}{4}(m-1).
\end{align*}

Hence $\delta(H^1_{\pi})\ge \frac{3}{4}(m-1)$ and then Theorem~\ref{thm:bipartite spread} implies a $(C_{\ref{thm:bipartite spread}}/(m-1))$-spread distribution on perfect matchings of $H^1_{\pi}$.
Sample a perfect matching $M_{\pi}^{1}$ from the $(C_{\ref{thm:bipartite spread}}/(m-1))$-spread distribution (on perfect matchings) of $H^1_{\pi}$.
For each $D_{i}\in B$, let $v^{1}_{i}$ be the vertex of $A$ that is matched to $D_{i}$ in $M_{\pi}^{1}$.
Then let $D_{i}^{1}=D_{i}\cup \{v^{1}_{i}\}$ for each $D_{i}\in B$. This completes the first round.

Now we do the second round of the distribution by defining another balanced bipartite graph $H^2_{\pi}$ with two parts $A',B'$ as follows: $A':={\textstyle \bigcup_{D_{i}\in \mathcal{D}_{\pi}}^{}}D_{i}\cap V_{2}$, $B':=\{D_{i}^1: D_{i}\in B\}$ and add an edge between $v^{2}\in A'$ and $D_{i}^{1}\in B'$ if $\delta_{\{1\}}(G[D_{i}^{1}\cup \{v^{2}\}])+\delta_{\{2\}}(G[D_{i}^{1}\cup \{v^{2}\}])\ge(1+\frac{\varepsilon}{3})C$, i.e., 
\begin{align*}
\delta_{\{1\}}(G[D_{i}\cup \{v_{i}^{1}, v^{2}\}])+\delta_{\{2\}}(G[D_{i}\cup \{v_{i}^{1}, v^{2}\}])\ge(1+\tfrac{\varepsilon}{3})C.
\end{align*}
Recall that in the definition of $H_{\pi}^{1}$, for every $D_{i}^{1}\in B'$, there are at least $(1-2\exp(-C\varepsilon^{2}/16))n$ vertices $v^{2}\in V_{2}$ such that $\delta_{\{1\}}(G[D_{i}^1\cup v^2])+\delta_{\{2\}}(G[D_{i}^1\cup v^2])\ge(1+\frac{\varepsilon}{3})C$.
So we have 
\begin{align*}
d_{H_{\pi}^2}(D_{i}^1)\ge (m-1)-2\exp(-C\varepsilon^{2}/16)n\ge \tfrac{3}{4}(m-1).
\end{align*}
For every $v^2\in A'$,  by Lemma~\ref{lem:auxiliary bipartite degree 1}, there are at least $(1-1/C)m'-|\mathcal{D}_{\pi}|$ clusters $D_{i}$ such that $\delta_{\{1\}}(G[D_{i}^1\cup \{v^2\}])+\delta_{\{2\}}(G[D_{i}^1\cup \{v^2\}])\ge(1+\tfrac{\varepsilon}{3})C$.
So we have 
\begin{align*}
d_{H_{\pi}^2}(v^2)\ge (1-1/C)m'-|\mathcal{D}_{\pi}|=(1-1/C)m'-(m'-m)\ge \tfrac{3}{4}(m-1).
\end{align*}
Hence $\delta(H^2_{\pi})\ge \frac{3}{4}(m-1)$. By Theorem~\ref{thm:bipartite spread}, there is a $(C_{\ref{thm:bipartite spread}}/(m-1))$-spread distribution on perfect matchings of $H_{\pi}^{2}$.

In the same way, we distribute vertices of $A'$ by defining the random sets $U_{i}$ as follows.
Sample $M_{\pi}^{2}$ from the $(C_{\ref{thm:bipartite spread}}/(m-1))$-spread distribution on perfect matchings of $H_{\pi}^{2}$.
For each $D_{i}^{1}\in B'$, let $v_i^{2}\in A'$ be the unique vertex that is matched to $D_{i}^{1}$ in $M_{\pi}^{2}$ and let $U_{i}:=D_{i}^{1}\cup \{v_i^{2}\}=D_{i}\cup \{v_i^{1}, v_i^{2}\}$ for each $i\in \{2,\dots , m\}$, whilst $U_{1}:=D_{1}$. 
Now it remains to check the resulting partition $\mathcal{U}=\{U_{1}, U_{2}, \dots, U_{m}\}$ satisfying the properties~\ref{size2}\ref{degree2}\ref{location2} in Lemma~\ref{lem:bipartite clustering}.
Since~\ref{size2} and~\ref{degree2} easily follow from our construction, we only verify~\ref{location2}.

Let $y_{1}, y_{2}, \dots ,y_{s}\in V(G)$, and $f: [s]\rightarrow [m]$ be given as in~\ref{location2}.
Let $T_{i}$ be the event that $y_{i}\in U_{f(i)}$, and $T_{i}^{1}$ be the event that $y_{i}\in D_{f(i)}$ whilst $T_{i}^{2}$ be the event that $y_{i}\in \{v_{f(i)}^{1}, v_{f(i)}^{2}\}$. Note that $T_i=T_i^1\cup T_i^2$ and the event $T_{i}^{1}$ is determined only by the permutation $\pi=(\pi_{1},\pi_{2})$, while the event $T_{i}^{2}$ happens if $\{y_i,D_{f(i)}\}$ is an edge of $M_\pi^1$ or $\{y_i,D_{f(i)}^1\}$ is an edge of $M_\pi^2$.
Let $S=\{i\in [s]:T_i^1~\text{happens}\}$, $S_1=S\cap \{i\in [s]:y_i\in V_1\}$, $S_2=S\setminus S_1$, $S_3=([s]\setminus S)\cap \{i\in [s]:y_i\in V_1\}$ and $S_4=([s]\setminus S)\cap \{i\in [s]:y_i\in V_2\}$.
By Stirling's formula we have
\begin{align*}
\mathbb{P}\left[\bigcap_{i\in S}^{} T_{i}^{1}\right]
\le \frac{C^{2|S|}\cdot (n-|S_{1}|)!\cdot (n-|S_{2}|)!}{n!\cdot n!}/\mathbb{P}[E_{1}\cap E_{2}\cap E_{3}]
\le \frac{100}{97}\left(\frac{C^{2}e}{n}\right)^{|S|}.
\end{align*}
Denote by $T_{i,1}^{2}$ the event $y_{i}=v_{f(i)}^{1}$ and by $T_{i,2}^{2}$ the event $y_{i}=v_{f(i)}^{2}$.
Note that $T_{i}^{2}=T_{i,1}^{2}\cup T_{i,2}^{2}$.
For every $\pi=(\pi_{1}, \pi_{2})\in E_{1}\cap E_{2}\cap E_{3}$, 
by the distribution on perfect matchings of bipartite graphs $H_{\pi}^{1}$ and $H_{\pi}^{2}$, we have
\begin{align*}
\mathbb{P}\left[\bigcap_{i\in [s]\setminus S}^{} T_{i}^{2} \Big| (\pi_{1}, \pi_{2} )\right]
=&\mathbb{P}\left[\bigcap_{i\in S_3}^{} T_{i,1}^{2}\Big| (\pi_{1}, \pi_{2})\right]\cdot \mathbb{P}\left[\bigcap_{i\in S_4}^{} T_{i,2}^{2}\Big| \bigcap_{i\in S_3}^{} T_{i,1}^{2}, (\pi_{1}, \pi_{2})\right]\\
\le &\left(\frac{C_{\ref{thm:bipartite spread}}}{m-1}\right)^{|S_{3}|}\cdot \left(\frac{C_{\ref{thm:bipartite spread}}}{m-1}\right)^{|S_{4}|}\\
\le & \left(\frac{C_{\ref{thm:bipartite spread}}}{m-1}\right)^{s-|S|}.
\end{align*}
Therefore, 
\begin{align*}
\mathbb{P}\left[\bigcap_{i\in S}T_{i}^{1}\cap \bigcap_{i\in [s]\setminus S}^{} T_{i}^{2}\right]
=\mathbb{P}\left[\bigcap_{i\in S}^{} T_{i}^{1}\right]\cdot \mathbb{P}\left[ \bigcap_{i\in [s]\setminus S}^{} T_{i}^{2} \Big| \bigcap_{i\in S}^{} T_{i}^{1}\right]
\le \frac{100}{97}\left(\frac{C^{2}e}{n}\right)^{|S|} \cdot\left(\frac{C_{\ref{thm:bipartite spread}}}{m-1}\right)^{s-|S|}
\le \left(\frac{C'}{2n}\right)^{s}.
\end{align*}
Finally, by the union bound over the $2^{s}$ choices of $S\< [s]$, we have
\begin{align*}
&\mathbb{P}\left[y_{i}\in U_{f(i)} \ \text{for every} \ i\in [s]\right]=\mathbb{P}\left[\bigcap_{i\in [s]} T_{i}\right]
=\sum_{S\< [s]}\mathbb{P}\left[\bigcap_{i\in S}^{} T_{i}^{1}\cap \bigcap_{i\in [s]\setminus S}^{} T_{i}^{2}\right]
\le \left(\frac{C'}{n}\right)^{s}.
\end{align*}
\end{proof}

\section{Concluding remarks}
In an earlier version of this paper~\cite{Han-Hu-Yang}, we derived the robustness of multipartite {H}ajnal--{S}zemer{\'e}di theorem. 
Indeed, by the random clustering method, we can obtain the robustness of existence of $F$-factors in multipartite graphs and its transversal version, where $F$ is a strictly 1-balanced graph. Before stating the results, we give some definitions.

Let $G$ be an $r$-partite graph with vertex classes $V_1,\ldots,V_r$.
We say that $G$ is \textit{balanced} if $|V_i|=|V_j|$ for any $1\le i < j\le r$.
Write $G[V_i,V_j]$ for the induced bipartite subgraph on vertex classes $V_i$ and $V_j$.
Define $\delta^\ast(G)$ to be $\min_{1\le i<j\le r} \delta(G[V_i,V_j])$.
Let $\mathbf{H} =\{H_{1}, H_{2}, \dots ,H_{m}\}$ be any graph system on $rn$ vertices where each $H_i$ is a balanced $r$-partite graph with each vertex class of size $n$, and we define $\de^*(\mathbf{H})=\min_{i\in [m]}\de^*(H_i)$.

\begin{defn}(multipartite Dirac threshold)
Let $F$ be a graph on $s$ vertices with $\chi(F)=r$.
By $\delta^{*}_{F}$ we denote the infimum real number $\delta$ such that for all $\alpha> 0$ and sufficiently large $n$ the following holds.
Let $G$ be any balanced $r$-partite graph with each vertex class of size $n$ with $n\in s\mathbb{N}$.
If $\delta^{*}(G)\ge (\delta+\alpha)n$, then $G$ contains an $F$-factor.
\end{defn}

\begin{theorem}\label{thm:multipartite}
For every $\alpha>0$ and $s,r\in\N$, there exists $C=C(s,r,\alpha)>0$ such that the following holds. Let $F$ be a strictly $1$-balanced graph on $s$ vertices with $\chi(F)=r$, and $G$ be a balanced $r$-partite graph with each vertex class of size $n$ such that $\delta^{*}(G)\ge (\delta^{*}_{F}+\alpha)n$. If $n\in s\mathbb{N}$ and $p\ge Cn^{-1/d_1(F)}(\log n)^{1/|E(F)|}$, then w.h.p.~$G(p)$ contains an $F$-factor.
\end{theorem}

\begin{defn}(multipartite Dirac threshold for graph systems)
Let $F$ be a graph on $s$ vertices with $t$ edges and $\chi(F)=r$.
By $\delta_{F}^{T*}$ we denote the infimum real number $\delta$ such that for all $\alpha> 0$ and sufficiently large $n$ the following holds.
Let $\mathbf{H} =\{H_{1}, H_{2}, \dots ,H_{trn}\}$ be any graph system on $srn$ vertices where each $H_i$ is a balanced $r$-partite graph with each vertex class of size $sn$.
If $\delta^{*}(\mathbf{H})\ge (\delta+\alpha)sn$, then $\mathbf{H}$ contains a transversal $F$-factor.
\end{defn}

\begin{theorem}\label{thm:transversal multipartite}
For every $\alpha>0$ and $r,s,t\in\N$, there exists $C=C(r,s,t,\alpha)>0$ such that the following holds. 
Let $F$ be a strictly 1-balanced graph on $s$ vertices with $t$ edges and $\chi(F)=r$, and $\mathbf{H} =\{H_{1}, H_{2}, \dots ,H_{trn}\}$ be any graph system on $srn$ vertices with $\delta^{*}(\mathbf{H})\ge (\delta_{F}^{T*}+\alpha)sn$. If $p\ge Cn^{-1/d_1(F)-1}(\log n)^{1/t}$, then w.h.p.~$\mathbf{H}(p) :=\{H_{1}(p), H_{2}(p), \dots ,H_{trn}(p)\}$ contains a transversal $F$-factor.
\end{theorem}

We do not give the proofs of Theorems~\ref{thm:multipartite} and~\ref{thm:transversal multipartite} here, as they are similar to the proof of Theorem~\ref{thm: trans 1-bal} and will be covered in the third author's master thesis.

The following result is a corollary of Theorem~\ref{thm: trans 1-bal}.

\begin{cor}\label{cor:trans 1-bal}
For every $\alpha>0$ and $k,s,t\in\N$, there exists $C=C(k,s,t,\alpha)>0$ such that the following holds for every $d\in[k-1]$. 
Let $F$ be a strictly $1$-balanced $k$-graph on $s$ vertices with $t$ edges, and $H$ be any $k$-graph on $sn$ vertices with $\delta_{d}(H)\ge (\delta^T_{F,d}+\alpha)\binom{sn-d}{k-d}$.
If $p\ge Cn^{-1/d_1(F)-1}(\log n)^{1/t}$, then w.h.p.~$\mathbf{H} :=\{H_{1}, H_{2}, \dots ,H_{tn}\}$ contains a transversal $F$-factor where each $H_i$ is independently distributed as $H(p)$.
\end{cor}

The $k$-graph system in Corollary~\ref{cor:trans 1-bal} consists of independent random subhypergraphs of a given $k$-graph $H$. What if we intersect each hypergraph of a given $k$-graph system with the same binomial random $k$-graph $\mathcal{G}^{(k)}(n, p)$?

\begin{question}\label{question}
For every $\alpha>0$ and $k,s,t\in\N$, there exists $C=C(k,s,t,\alpha)>0$ such that the following holds for every $d\in[k-1]$. 
Let $F$ be a strictly $1$-balanced $k$-graph on $s$ vertices with $t$ edges, and $\mathbf{H} =\{H_{1}, H_{2}, \dots ,H_{tn}\}$ be any $k$-graph system on $sn$ vertices with $\delta_{d}(\mathbf{H})\ge (\delta^T_{F,d}+\alpha)\binom{sn-d}{k-d} $. 
Suppose that $p\ge Cn^{-1/d_1(F)}(\log n)^{1/t}$.
Is there a transversal $F$-factor in $\mathbf{H}\cap \mathcal{G}^{(k)}(n, p)$ with high probability?
\end{question}

Note that we put the bound $p\ge Cn^{-1/d_1(F)}(\log n)^{1/t}$, which is necessary up to the value of $C$ even if all $H_i$ are complete graphs. 
Moreover, we can show that Question~\ref{question} for $F=K_r$ has a positive answer following the method of Ferber, Han and Mao~\cite{Ferber}, which shows some evidence for the general case.

\bibliographystyle{abbrv}
\normalem
\bibliography{ref}

\begin{thebibliography}{10}

\bibitem{MR4125343}
R.~Aharoni, M.~DeVos, S.~Gonz\'{a}lez Hermosillo de~la Maza, A.~Montejano, and
  R.~\v{S}\'{a}mal.
\newblock A rainbow version of {M}antel's theorem.
\newblock {\em Adv. Comb.}, pages Paper No. 2, 12pp, 2020.

\bibitem{Aharoni}
R.~Aharoni, A.~Georgakopoulos, and P.~Spr{\"u}ssel.
\newblock Perfect matching in $r$-partite $r$-graphs.
\newblock {\em Eur. J. Comb.}, 30:39--42, 2009.

\bibitem{MR3628907}
R.~Aharoni and D.~Howard.
\newblock A rainbow {$r$}-partite version of the {E}rd{\H{o}}s-{K}o-{R}ado
  theorem.
\newblock {\em Combin. Probab. Comput.}, 26(3):321--337, 2017.

\bibitem{Allen}
P.~Allen, J.~B\"ottcher, J.~Corsten, E.~Davies, M.~Jenssen, P.~Morris,
  B.~Roberts, and J.~Skokan.
\newblock A robust {C}orr\'adi-{H}ajnal theorem.
\newblock {\em Random Struct. Algorithms}, 65(1):61--130, 2024.

\bibitem{Anastos-Chakraborti}
M.~Anastos and D.~Chakraborti.
\newblock Robust {H}amiltonicity in families of {D}irac graphs.
\newblock {\em arXiv preprint arXiv:2309.12607}, 2023.

\bibitem{Bastide}
P.~Bastide, C.~Legrand-Duchesne, and A.~M{\"u}yesser.
\newblock Random embeddings of bounded degree trees with optimal spread.
\newblock {\em arXiv preprint arXiv:2409.06640}, 2024.

\bibitem{MR4394673}
P.~Bradshaw.
\newblock Transversals and bipancyclicity in bipartite graph families.
\newblock {\em Electron. J. Combin.}, 28(4):Paper No. 4.25, 20pp, 2021.

\bibitem{rainbow-bandwidth}
D.~Chakraborti, S.~Im, J.~Kim, and H.~Liu.
\newblock A bandwidth theorem for graph transversals.
\newblock {\em arXiv preprint arXiv:2302.09637}, 2023.

\bibitem{CHWW1}
Y.~Cheng, J.~Han, B.~Wang, and G.~Wang.
\newblock Rainbow spanning structures in graph and hypergraph systems.
\newblock {\em Forum Math. Sigma}, 11:Paper No. e95, 20, 2023.

\bibitem{CHWW2}
Y.~Cheng, J.~Han, B.~Wang, G.~Wang, and D.~Yang.
\newblock Transversal {H}amilton cycle in hypergraph systems.
\newblock {\em SIAM J. Discrete Math.}, 39(1):55--74, 2025.

\bibitem{MR4287703}
Y.~Cheng, G.~Wang, and Y.~Zhao.
\newblock Rainbow pancyclicity in graph systems.
\newblock {\em Electron. J. Combin.}, 28(3):Paper No. 3.24, 9pp, 2021.

\bibitem{Corradi}
K.~Corr{\'a}di and A.~Hajnal.
\newblock On the maximal number of independent circuits in a graph.
\newblock {\em Acta Math. Acad. Sci. Hungar.}, 14:423--439, 1963.

\bibitem{Dirac}
G.~A. Dirac.
\newblock Some theorems on abstract graphs.
\newblock {\em Proc. Lond. Math. Soc.}, 3(1):69--81, 1952.

\bibitem{Ferber}
A.~Ferber, J.~Han, and D.~Mao.
\newblock Dirac-type problem of rainbow matchings and {H}amilton cycles in
  random graphs.
\newblock {\em arXiv preprint arXiv:2211.05477}, 2022.

\bibitem{Frankston}
K.~Frankston, J.~Kahn, B.~Narayanan, and J.~Park.
\newblock Thresholds versus fractional expectation thresholds.
\newblock {\em Ann. of Math.}, 194(2):475--495, 2021.

\bibitem{Gupta}
P.~Gupta, F.~Hamann, A.~M{\"u}yesser, O.~Parczyk, and A.~Sgueglia.
\newblock A general approach to transversal versions of {D}irac-type theorems.
\newblock {\em Bull. Lond. Math. Soc}, 55(6):2817--2839, 2023.

\bibitem{HSz}
A.~Hajnal and E.~Szemer{\'e}di.
\newblock Proof of a conjecture of {P}. {E}rd{\H{o}}s.
\newblock {\em Comb. Theory Appl.}, 2(4):601--623, 1970.

\bibitem{HPS}
H.~H{\`a}n, Y.~Person, and M.~Schacht.
\newblock On perfect matchings in uniform hypergraphs with large minimum vertx
  degree.
\newblock {\em SIAM J. Discret. Math.}, 23(2):732--748, 2009.

\bibitem{Han-Hu-Yang}
J.~Han, J.~Hu, and D.~Yang.
\newblock A robust version of the multipartite {H}ajnal--{S}zemer{\'e}di
  theorem.
\newblock {\em arXiv preprint arXiv:2311.00950v1}, 2023.

\bibitem{Janson}
S.~Janson, T.~{\L}uczak, and A.~Rucinski.
\newblock {\em Random graphs}.
\newblock Wiley-Interscience Series in Discrete Mathematics and Optimization.
  Wiley-Interscience, New York, 2000.

\bibitem{Johansson}
A.~Johansson, J.~Kahn, and V.~Vu.
\newblock Factors in random graphs.
\newblock {\em Random Struct. Algorithms}, 33(1):1--28, 2008.

\bibitem{MR4171383}
F.~Joos and J.~Kim.
\newblock On a rainbow version of {D}irac's theorem.
\newblock {\em Bull. Lond. Math. Soc.}, 52(3):498--504, 2020.

\bibitem{Joos-Lang}
F.~Joos, R.~Lang, and N.~Sanhueza-Matamala.
\newblock Robust {H}amiltonicity.
\newblock {\em arXiv preprint arXiv:2312.15262}, 2023.

\bibitem{Kahn-Kalai}
J.~Kahn and G.~Kalai.
\newblock Thresholds and expectation thresholds.
\newblock {\em Combin. Probab. Comput.}, 16(3):495--502, 2007.

\bibitem{Kang-Kelly}
D.~Y. Kang, T.~Kelly, D.~K{\"u}hn, D.~Osthus, and V.~Pfenninger.
\newblock Perfect matchings in random sparsifications of {D}irac hypergraphs.
\newblock {\em Combinatorica}, 44:1233--1266, 2024.

\bibitem{KMP}
T.~Kelly, A.~M\"uyesser, and A.~Pokrovskiy.
\newblock Optimal spread for spanning subgraphs of {D}irac hypergraphs.
\newblock {\em J. Combin. Theory Ser. B}, 169:507--541, 2024.

\bibitem{MR3180741}
M.~Krivelevich, C.~Lee, and B.~Sudakov.
\newblock Robust {H}amiltonicity of {D}irac graphs.
\newblock {\em Trans. Amer. Math. Soc.}, 366(6):3095--3130, 2014.

\bibitem{DKDO}
D.~K{\"u}hn and D.~Osthus.
\newblock Matchings in hypergraphs of large minimum degree.
\newblock {\em J. Graph Theory}, 51(4):269--280, 2006.

\bibitem{MR2506388}
D.~K\"{u}hn and D.~Osthus.
\newblock The minimum degree threshold for perfect graph packings.
\newblock {\em Combinatorica}, 29(1):65--107, 2009.

\bibitem{Liebenau}
A.~Liebenau and N.~Wormald.
\newblock Asymptotic enumeration of graphs by degree sequence, and the degree
  sequence of a random graph.
\newblock {\em J. Eur. Math. Soc.}, 26(1):1--40, 2023.

\bibitem{MR4451911}
H.~Lu, Y.~Wang, and X.~Yu.
\newblock Rainbow perfect matchings for 4-uniform hypergraphs.
\newblock {\em SIAM J. Discrete Math.}, 36(3):1645--1662, 2022.

\bibitem{MR4523452}
H.~Lu, Y.~Wang, and X.~Yu.
\newblock A better bound on the size of rainbow matchings.
\newblock {\em J. Combin. Theory Ser. A}, 195(105700), 2023.

\bibitem{MR4275007}
H.~Lu, X.~Yu, and X.~Yuan.
\newblock Rainbow matchings for 3-uniform hypergraphs.
\newblock {\em J. Combin. Theory Ser. A}, 183(105489), 2021.

\bibitem{LLL}
L.~Lu, P.~Li, and X.~Li.
\newblock Rainbow structures in a collection of graphswith degree conditions.
\newblock {\em J. Graph Theory}, 104(2):341--359, 2023.

\bibitem{Michael-Bruce}
M.~Molloy and B.~Reed.
\newblock {\em Graph colouring and the probabilistic method}, volume~23 of {\em
  Algorithms and Combinatorics}.
\newblock Springer-Verlag, Berlin, 2002.

\bibitem{MR4451150}
R.~Montgomery, A.~M\"{u}yesser, and Y.~Pehova.
\newblock Transversal factors and spanning trees.
\newblock {\em Adv. Comb.}, pages Paper No. 3, 25pp, 2022.

\bibitem{Park-Pham}
J.~Park and H.~T. Pham.
\newblock A proof of the {K}ahn-{K}alai conjecture.
\newblock {\em J. Amer. Math. Soc.}, 37(1):235--243, 2024.

\bibitem{Pham2}
H.~T. Pham, A.~Sah, M.~Sawhney, and M.~Simkin.
\newblock A toolkit for robust thresholds.
\newblock {\em arXiv preprint arXiv:2210.03064v3}, 2022.

\bibitem{Pikhurko}
O.~Pikhurko.
\newblock Perfect matching and {$K_{4}^{3}$}-tilings in hypergraphs of large
  codegree.
\newblock {\em Graphs and Combinatorics}, 24:391--404, 2008.

\bibitem{Riordan}
O.~Riordan.
\newblock Random cliques in random graphs and sharp thresholds for
  {$F$}-factors.
\newblock {\em Random Struct. Algorithms}, 61(4):619--637, 2022.

\bibitem{Rodl}
V.~R{\"o}dl and A.~Ruci{\'n}ski.
\newblock Dirac-type questions for hypergraphs -- {A} survey (or more problems
  for {E}ndre to solve).
\newblock {\em Bolyai Soc. Math. Stud.}, 21:561--590, 2010.

\bibitem{Rucinski1992}
A.~Ruci\'nski.
\newblock Matching and covering the vertices of a random graph by copies of a
  given graph.
\newblock {\em Discrete. Math.}, 105:185--197, 1992.

\bibitem{Sudakov-survey}
B.~Sudakov.
\newblock Robustness of graph properties.
\newblock {\em Surveys in combinatorics 2017, London Math. Soc. Lecture Note
  Ser.}, Cambridge Univ. Press, Cambridge, 440:372--408, 2017.

\bibitem{Sudakov2008}
B.~Sudakov and V.~H. Vu.
\newblock Local resilience of graphs.
\newblock {\em Random Struct. Algorithms}, 33(4):409--433, 2008.

\bibitem{SWW}
W.~Sun, G.~Wang, and L.~Wei.
\newblock Transversal structures in graph systems: A survey.
\newblock {\em arXiv preprint arXiv.2412.01121}, 2024.

\bibitem{Talagrand}
M.~Talagrand.
\newblock Are many small sets explicitly small?
\newblock In {\em Proceedings of the forty-second ACM symposium on Theory of
  computing}, pages 13--36, 2010.

\end{thebibliography}

\begin{appendices}
\section{Proof of Theorem~\ref{thm:multipartite spread}}\label{app1}
We need the following concentration result in~\cite{Pikhurko}.
\begin{lemma}\label{cor1}
Let $G$ be an arbitrary subgraph of $K_{n,n}$ and let $h$ satisfy $h\geq 4np_0$, where $p_0=4\exp(-h^2/(256n))$. Let $M$ be a perfect matching of $K_{n,n}$ chosen uniformly at random. Then 
\[
\mathbb{P}\left[\big{|}|M\cap E(G)|-\tfrac{e(G)}{n}\big{|}\ge h\right]\le p_0.
\]  
\end{lemma}

\begin{proof}[Proof of Theorem~\ref{thm:multipartite spread} by using Theorem~\ref{thm:bipartite perfect matching}]
Given $\eps>0$ and $k\ge 3$, we choose $1/n\ll 1/C\ll \eps,1/k$. Let $L\in \binom{[k]}{\ell}$ and
$H$ be a $k$-partite $k$-graph with parts $V_{1}\cup \dots \cup V_{k}=V(H)$ such that $|V_{i}|=n$ for each $i\in [k]$ and 
\[\delta_{L}(H)n^{\ell}+\delta_{[k]\setminus L}(H)n^{k-\ell}\ge (1+\varepsilon )n^{k}.\]
Assume without loss of generality that $L=[\ell]$ and $\ell\ge 2$. We will inductively produce perfect matchings $M_{1}, \dots , M_{\ell -1}$ where $M_{i}$ is a perfect matching chosen uniformly at random in the complete bipartite graph $K[V_{i}, V_{i+1}]$.
Given $i\in [\ell -1]$ and $M_{1}, \dots , M_{i-1}$, let $T_{i}$ be a perfect matching in the complete $i$-partite $i$-graph $K[V_{1}, \dots, V_{i}]$ obtained by ``gluing" the edges of $M_{1}, \dots , M_{i-1}$.
Formally, $T_{1}$ is the set of vertices in $V_{1}$, which is trivial.
For $i\ge2$,
\begin{align*}
T_{i}=\{\{x_{1},\dots,x_{i}\}: \forall j\in [i-1], \{x_{j},x_{j+1}\}\in M_{j}\}.
\end{align*}

We call a set $X\< V_{i+2}\cup\dots\cup V_{k}$ \textit{legal} if $|X\cap V_{j}|=1$ for every $j\in \{i+2, \dots , k\}$. For any legal set $X\< V_{i+2}\cup\dots\cup V_{k}$, let $F_{X}\< K[V_{i}, V_{i+1}]$ be the set of pairs $\{x_{i}, x_{i+1}\}$ with $x_{i}\in V_{i}$ and $x_{i+1}\in V_{i+1}$ such that $D\cup \{x_{i+1}\} \cup X\in E(H)$, where $D$ is an element of $T_{i}$ that contains $x_{i}$.
Suppose we have already constructed $M_{1}, \dots ,M_{i-1}$.
Let $M_{i}$ be a perfect matching of $K[V_{i},V_{i+1}]$ chosen uniformly at random.
By Lemma~\ref{cor1} and the union bound,
we know that with probability $1-o(\tfrac{1}{k})$, $M_{i}$ satisfies that
\begin{align}\label{eq1}
|M_{i}\cap F_{X}|=\frac{|F_{X}|}{n}\pm \lambda \ \text{for every}\ X\in K[V_{i+2},\dots , V_{k}],
\end{align}
where $\lambda=\sqrt{257kn\log n}$.
We call a perfect matching $M_{i}$ satisfying this property to be ``good". Fix a good $M_i$ uniformly at random and repeat this for $M_{i+1}$ (until $i=\ell-1$).

Hence for every $i\in [\ell-1]$ there is a $(C/n)$-vertex-spread distribution on embeddings of ``good" perfect matchings $M_i$ in $K[V_{i}, V_{i+1}]$. Next we only consider good perfect matchings, and conditioning on this,
we claim that for every $i\in [\ell -1]$ and every legal set $X\< V_{i+2}\cup\dots\cup V_{k}$, it holds that
\begin{align}\label{eq2}
|N(X)\cap T_{i+1}|=\frac{|N(X)|}{n^{i}}\pm i\lambda.
\end{align}
We prove this by induction on $i$, and the basic case $i=1$ just follows from \eqref{eq1}.
Suppose that $2\le i\le \ell-1$ and the inequality \eqref{eq2} holds for $i-1$.
For every legal $X\< V_{i+2}\cup\dots\cup V_{k}$, we have
\begin{align*}
|N(X)\cap T_{i+1}|&=|F_{X}\cap M_{i}|=\frac{1}{n}|F_{X}|\pm \lambda\\
&=\frac{1}{n}\sum_{x\in V_{i+1}}|N(X\cup\{x\})\cap T_{i}|\pm \lambda\\
&=\frac{1}{n}\sum_{x\in V_{i+1}}\left(\frac{|N(X\cup\{x\})|}{n^{i-1}}\pm(i-1)\lambda\right)\pm \lambda\\
&=\frac{|N(X)|}{n^{i}}\pm i\lambda,
\end{align*}
where the second equality follows as $M_{i}$ is good and the penultimate equality follows by the induction hypothesis.
This tells that for every legal $X\< V_{\ell+1}\cup\dots\cup V_{k}$ we have
\begin{align}\label{eq3}
|N(X)\cap T_{\ell }|\ge \frac{|N(X)|}{n^{\ell-1}}-(\ell-1)\lambda\ge \frac{\delta_{[k]\setminus L}(H)}{n^{\ell-1}}-(\ell-1)\lambda.
\end{align}

Similarly, we inductively build perfect matchings $M_{j}'\< K[V_{j-1}, V_{j}]$ for $j=k, k-1, \dots , \ell+2$ chosen uniformly at random.
Similar arguments also, with probability $1-o(1)$, yield a perfect matching $T_{k-\ell}'$ in $K[V_{\ell+1},\dots , V_{k}]$ such that every legal $\ell$-set $Y\<V_{1}\cup\dots\cup V_{\ell}$ satisfies
\begin{align}\label{eq4}
|N(Y)\cap T_{k-\ell }'|\ge \frac{\delta_{L}(H)}{n^{k-\ell-1}}-(k-\ell-1)\lambda.
\end{align}

Now consider an auxiliary balanced bipartite graph $G$ with parts $U_{1}=T_{\ell}$ and $U_{2}=T_{k-\ell}'$, where $X\in T_{\ell}$ and $Y\in T_{k-\ell}'$ are adjacent in $G$ whenever $X\cup Y\in E(H)$.
Combining \eqref{eq3} and \eqref{eq4}, we have
\begin{align*}
\delta_{\{1\}}(G)+\delta_{\{2\}}(G)\ge \frac{\delta_{[k]\setminus L}(H)}{n^{\ell-1}}+\frac{\delta_{L}(H)}{n^{k-\ell-1}}-(\ell-1)\lambda-(k-\ell-1)\lambda\ge (1+\eps/2)n.
\end{align*}
Then by Theorem~\ref{thm:bipartite perfect matching}, there is a $(C/n)$-vertex-spread distribution on embeddings of perfect matchings in $G$. Observe that every perfect matching $M_G$ in $G$, together with the perfect matchings $T_{\ell}$ and $T'_{k-\ell}$, actually forms a perfect matching in $H$, denoted as $(T_{\ell},M_G,T'_{k-\ell})$. This immediately gives an embedding $\psi$ of perfect matching $M$ into $H$ as follows:
\begin{enumerate}[label = (\roman{enumi})]
    \item assume $V(M)=[kn]$ and for each $p\in[n]$, $\{p,n+p,2n+p,\ldots, (k-1)n+p\}$ is an edge in $M$; 
    \item \label{ii} the restrictions $\psi:[n]\rightarrow V_1$ and $[(k-1)n+1,kn]\rightarrow V_k$ are  the bijections chosen uniformly and independently, and for the special case $k-\ell=1$, we just solely choose the former restriction;
    \item for each $p\in[n]$, $i\in[2,\ell]$ and $j\in[\ell+1,k-1]$, $\psi(\{p,n+p,2n+p,\ldots, (i-1)n+p\})$ is the unique element in $T_i$, while $\psi(\{(j-1)n+p,jn+p,\ldots, (k-1)n+p\})$ is the unique element in $T'_{k-j+1}$;
    \item for each $p\in[n]$, $\psi(\{p,n+p,2n+p,\ldots, (\ell-1)n+p\})$ and $\psi(\{\ell n+p,\ldots, (k-1)n+p\})$ are adjacent in $M_G$.
\end{enumerate}

We claim that such an embedding $\psi$ admits a $(C/n)$-vertex-spread distribution. That is, for every $s\in [kn]$ and every two distinct sequences $x_{1},\dots, x_{s}\< V(M)$ and $y_{1},\dots, y_{s}\< V(H)$, we have
\[\mathbb{P}\big{[}\psi(x_{i})=y_{i} \ \text{for all}\ i\in [s]\big{]}\le (C/n)^{s}.\]
Let $I_{i}$ be the index set of all vertices in $\{y_{1},\dots, y_{s}\}\cap V_{i}$, where $i\in[k]$ and $E_i:=\bigcap_{j\in I_i}\{\psi(x_{j})=y_{j}\}$. Then \[\mathbb{P}\big{[}\psi(x_{i})=y_{i} \ \text{for all}\ i\in [s]\big{]}=\mathbb{P}\big{[}\cap_{i\in[k]}E_i\big{]}.\]
For the case $k-\ell\ge 2$, observe that for $i\in[2,\ell]$ (or $i\in[\ell+1,k-1]$) the event $E_i$ is independently determined by the choice of the corresponding matching $M_{i-1}$ (resp. $M'_{k-i}$), while $E_1$ (or $E_k$) is determined by a  permutation uniformly chosen as in~\ref{ii}. Recall that any such $M_i$ admits a $(C/n)$-vertex-spread distribution, and therefore $\mathbb{P}[E_i]\le (C/n)^{|I_i|}$. This implies
\[\mathbb{P}\big{[}\cap_{i\in[k]}E_i\big{]}\le (\tfrac{C}{n})^{s}.
\]

For the special case $k-\ell=1$, we can see that the event $E_k=\bigcap_{j\in I_k}\{\psi(x_{j})=y_{j}\}$ is completely determined by the choice of the matching $M_{G}$. As $M_G$ also admits a $(C/n)$-vertex-spread distribution (conditional on the event that every $M_i$ is good for $i\in[\ell-1]$), we obtain
\begin{align*}
\mathbb{P}\big{[}\cap_{i\in[k]}E_i\big{]} 
=\mathbb{P}\big{[}E_k\big{|}\cap_{i\in[\ell]}E_i\big{]}\mathbb{P}\big{[}\cap_{i\in[\ell]}E_i\big{]}
\le (\tfrac{C}{n})^{|I_k|} (\tfrac{C}{n})^{s-|I_k|}= (\tfrac{C}{n})^{s}.
\end{align*}
\end{proof} 

\end{appendices}

\end{document}